\documentclass[12pt,letterpaper,reqno]{amsart}
\usepackage{amssymb}
\usepackage{pb-diagram}
\usepackage{pictexwd}
\usepackage{amsmath,amscd}
\usepackage{lamsarrow,pb-xy,pb-lams}
\usepackage{color}
\usepackage{epsf}

\newif\ifpdf
\ifpdf
 \usepackage[pdftex]{graphicx}
\else
 \usepackage[dvips]{graphicx}
\fi

\theoremstyle{plain}
\newtheorem{thm}{Theorem}[section]
\newtheorem*{main theorem}{Main Theorem}
\newtheorem{lemma}[thm]{Lemma}
\newtheorem{proposition}[thm]{Proposition}
\newtheorem{corollary}[thm]{Corollary}

\theoremstyle{definition}
\newtheorem{definition}[thm]{Definition}
\newtheorem{example}[thm]{Example}
\newtheorem{remark}[thm]{Remark}

\DeclareMathOperator{\Spin}{{\it Spin}}

\newcommand{\pardev}[2]{\frac{\partial #1}{\partial #2}}
\newcommand{\half}{\frac{1}{2}}

\newcommand{\R}{\mathbb R}
\newcommand{\Z}{\mathbb Z}

\newcommand{\N}{\mathbb N}
\newcommand{\F}{\mathbb F}

\newcommand{\co}[1]{\mathfrak{#1}^*}
\newcommand{\g}{\mathfrak g}
\newcommand{\Weil}[1]{W_\mathfrak{#1}}
\newcommand{\NWeil}[1]{\mathcal{W}_{\mathfrak{#1}}}
\newcommand{\sym}[1]{S(#1)}
\newcommand{\ext}[1]{\wedge(#1)}
\newcommand{\Cl}[1]{Cl(\mathfrak{#1})}
\newcommand{\U}[1]{\it{U}(\mathfrak{#1})}

\begin{document}

\title [\ ] {Superalgebraic interpretation of the quantization \\
       maps of Weil algebras}

\author[\ ]{Li Yu}

\thanks{Department of Mathematics, University of California, San Diego, La Jolla, California 92093-0112, USA. Email: lyu@math.ucsd.edu}
\date{May 1, 2005}

%TABLE OF CONTENTS, LISTS OF FIGURES & TABLES
%\tableofcontents

%Abstract
\begin{abstract}
   In \cite{AMM98}, A.Alekseev and E.Meinrenken construct
   an explicit $G$-differential space homomorphism $\mathcal{Q}$, called
   the quantization map, between the Weil algebra $\Weil{\g}=
   \sym{\co{\g}} \otimes \ext{\co{\g}}$ and $\NWeil{\g}=\U{\g} \otimes \Cl{\g}$
   (which they called the noncommutative Weil algebra) for any quadratic Lie algebra $\g$.
   They showed that $\mathcal{Q}$ induces an algebra isomorphism between
   the basic cohomology rings $H^{\ast}_{bas}(\Weil{g})$ and $H^{\ast}_{bas}(\NWeil{\g})$.
   In this paper, I will interpret the quantization map $\mathcal{Q}$
   as the super Duflo map between the symmetric algebra $S(\widetilde{T\g[1]})$
   and the universal enveloping algebra $U(\widetilde{T\g[1]})$ of a
   super Lie algebra $\widetilde{T\g[1]}$ which is canonically related to
   the quadratic Lie algebra $\g$.
   The basic cohomology rings $H^{\ast}_{bas}(\Weil{g})$ and $H^{\ast}_{bas}(\NWeil{g})$
   correspond exactly to $S(\widetilde{T\g[1]})^{inv}$ and $U(\widetilde{T\g[1]})$ respectively.
   So what they proved is equivalent to the fact that the Duflo map
   commutes with the adjoint action of the Lie algebra, and that the Duflo map
   is an algebra homomorphism when restricted to the space of invariants.
   In addition, I will explain how the diagrammatic analogue of the Duflo map
   introduced in \cite{BTT03} can be also made for the quantization map $\mathcal{Q}$.
\end{abstract}

\maketitle

%ACKNOWLEDGEMENTS
\textbf{Acknowledgements}:
First of all, I would like to thank my advisor J.Roberts and
professor N. Wallach for their support and
insightful advices.  I would also like to thank professor E. Meinrenken for pointing out some mistakes I wrote in the draft.
In addition, the argument in section 5.5 is due to Dylan Thurston.
I also understand that several people including E. Meinrenken and  P. Severa have already conjectured the main results of this paper and probably
have their way of proving them. However, no formal proof had yet been published. So I think it is worth writing it out. 

\large
\section{\textbf{Introduction}}

\normalsize

  The classical Weil algebra $W(\g)$ of a Lie algebra $\g$
  is introduced in the algebraic framework of equivariant geometry
  where people investigate the geometry and topology of group actions on
  smooth manifolds. As a vector space, $W(\g)$ is just
  $S(\co{\g}) \otimes \Lambda(\co{g})$. It has a well-known structure as a G-differential
  algebra and is very useful in equivariant De Rham theory.

  In 2000, A.Alekseev and E.Meinrenken introduced, for a Lie
  algebra with an ad-invariant metric (called a quadratic Lie algebra),
  the non-commutative version $\NWeil{\g}= U(g) \otimes Cl(g)$, which is
  also a G-differential algebra. They call it the noncommutative Weil algebra
  of $\g$. They then constructed a map
  $\mathcal{Q}: \Weil{\g} \longrightarrow \NWeil{\g}$,
  called the quantization map, between these two algebras which has
  three main properties:
  \begin{itemize}
  \item it is an isomorphism of vector spaces
  \item it is an isomorphism of G-differential spaces
  \item it is not a map of algebras, but it does induce an algebra isomorphism
        between the basic cohomologies of the two algebras.
  \end{itemize}
  However, their definition of $\mathcal{Q}$ and their proof are both quite complicated,
  and it is difficult to understand the geometric meaning of their formulae.

  Their theorem resembles another theorem, the Duflo isomorphism
  theorem in Lie theory. M.Duflo in \cite{Duf77} constructed, for any Lie algebra $\g$, a
  map $\Upsilon : S(\g) \longrightarrow U(\g)$ which has the
  properties :
  \begin{itemize}
   \item it is an isomorphism of vector spaces
   \item it is a map of $\g$-modules
   \item it is not a map of algebras, but induces an algebra isomorphism
         between the spaces of invariants.
  \end{itemize}
  This theorem is highly non-trivial although for semisimple Lie algebras, it follows from the
  Weyl character formula.

  In 2003 Dror Bar-Natan, Le Thang and Dylan Thurston proved in \cite{BTT03} yet another
  similar result, the "wheeling theorem". They constructed a map
  $\Phi: \mathcal{B} \longrightarrow \mathcal{A}$
  between certain spaces of diagrams appearing in Vassiliev theory,
  and proved that $\Phi$ is an algebra isomorphism with respect to
  some natural algebraic structures on these diagrams.

  Finally, there is Kontsevich's famous work \cite{Ko97} on deformation quantization
  of Poisson Manifolds. In his viewpoint, the Duflo map of a Lie algebra $\g$ comes from
  the deformation of the Lie-Poisson structure of $\g^*$. In addition, his proof can be
  generalized to the case of super Lie algebras.

  In this paper I will study the relationships between these
  theorems. The main result is:
  the Alekseev-Meinrenken quantization map $\mathcal{Q}$ can be identified with the
  Duflo map for a certain Lie superalgebra $\widetilde{Tg[1]}$, which is the
  Lie algebra of a central extension of the odd tangent bundle of G. Many properties of
  $\mathcal{Q}$ can therefore
  be proved by using Kontsevich's proof of the super-Duflo theorem, or by
  a diagrammatic proof using wheeling.

  Here are some logical relations between these theorems:
  \begin{itemize}
   \item Kontsevich's deformation quantization $\Rightarrow$ Duflo theorem
         (for any Lie superalgebras)
   \item Wheeling theorem $\Rightarrow$ Duflo theorem for quadratic Lie (super)algebras
   \item SuperDuflo $\Rightarrow$ Alekseev-Meinrenken theorem (this
   is the main theorem in this thesis)
  \end{itemize}

  The supergeometric interpretation of the non-commutative Weil algebra
  ought to be useful for future applications in geometry.

  \begin{remark}
   The noncommutative Weil algebra can so far only be defined for
   quadratic Lie algebras. It is unclear whether any of the
   results of Alekseev and Meinrenken can be extended to the
   non-quadratic case.
  \end{remark}

%%%%%%%%%%%%%%%%%%%%%%%%%%%%%%%%%%%%%%%%%%%%%%%%%%%%%%%%%%
%%%%%%%%% Plan of Paper %%%%%%%%%%%%%%%%%%%%%%%%%%%%%%%%%%
%%%%%%%%%%%%%%%%%%%%%%%%%%%%%%%%%%%%%%%%%%%%%%%%%%%%%%%%%%

\subsection{Plan of Paper}

I first present some standard introductory materials to make this
thesis as self-contained as possible in chapter 2. I will review
some basic facts about $G$-differential algebras and the Weil
algebra $\Weil{\g}$ of a Lie algebra $\g$ which are discussed
extensively in \cite{GS99}. A canonical $G$-differential structure
and its slightly varied form are defined on $\Weil{\g}$. In
addition I will define the $G$-differential structure on Clifford
algebra $\Cl{\g}$ and the \textit{noncommutative Weil algebra},
which is introduced in \cite{AM03} by A.Alekseev and E.Meinrenken.
Then I introduce the Duflo map for a Lie (super)algebra and the
quantization map of the Weil algebra $\Weil{\g}$.

In chapter 3, I will build a critical connection between the
$G$-differential structure on the Weil algebra $\Weil{\g}$ and the
Lie super algebra structure on $\widetilde{T\g[1]}$, which can
help us to understand the quantization map by the theory of Lie
algebras. Then I state the main theorem of this thesis.

In chapter 4, I will present a proof of the main theorem. To do
that, I have to first discuss spin representations for Clifford
algebras and its generalization to a super algebra. Then I
discuss the factorization of the spin representation constructed
by A.Alekseev and E.Meinrenken in \cite{AM02}. The proof of the
main theorem is put at the end of the chapter 4.

In Chapter 5, I will introduce Jacobi diagrams and diagrammatic
representation of tensors in the category of  Lie (super)
algebras. People can find the standard exposition of these in
\cite{Bar95} and \cite{BTT03}. In Chapter 6, I will explain how
the method of using Jacobi diagrams to prove the Duflo isomorphism
in \cite{BTT03} could be naturally extended to interpret the
quantization map.  \newline
 %introduction

\large

\section{\textbf{Definitions and Preliminary facts}}
\normalsize
%%%%%%%%%%%%%%%%%%%%%%%%%%%%%%%%%%%%%%%%%%%%%%%%%%%%%%%%
%%%%%%% G-differential algebras %%%%%%%%%%%%%%%%%%%%%%%%
%%%%%%%%%%%%%%%%%%%%%%%%%%%%%%%%%%%%%%%%%%%%%%%%%%%%%%%%

\subsection{G-differential algebras}
  Suppose $G$ is a Lie group with Lie algebra $\g$. Choose a basis
  $e_{1}, \ldots, e_{n}$ of $\g$ and let $e^{1}, \ldots, e^{n}$
  be the dual basis in $\co{\g}$. Let $\{f^{c}_{ab}\}$ be the structure constants
  defined by
  \[
      [e_{a},e_{b}]=f^{c}_{ab}e_{c}
  \]
  define $\hat{\g}$ to be the
  Lie super algebra
  \[
     \widehat{\g}:= \g_{-1} \oplus \g_{0} \oplus \g_{1}
  \]
  where $\g_{-1}$ is an $n$-dimensional vector space with basis
  $\iota_{1},\ldots, \iota_{n}$, $\g_{0}$ is an
  $n$-dimensional vector space with basis $L_{1}, \ldots, L_{n}$
  and $\g_{1}$ is a one-dimensional vector space with basis $d$. The
  super Lie bracket is defined in terms of this basis by
 \begin{align}
  [\iota_{a}, \iota_{b}] &= \iota_{a} \iota_{b}+ \iota_{b} \iota_{a}=0, \\
  [L_{a}, \iota_{b}]     &= L_{a}\iota_{b} - \iota_{b} L_{a}=f^{c}_{ab}\iota_{c}, \\
  [L_{a}, L_{b}]         &= L_{a}L_{b}-L_{b}L_{a}=f^{c}_{ab}L_{c}, \\
  [d, \iota_{a}]         &= d \iota_{a}+ \iota_{a}d = L_{a}, \\
  [d, L_{a}]              &= dL_{a} - L_{a}d=0 ,   \\
  [d, d]                 &= 2d^{2}=0.
 \end{align}

 Recall that a Lie super algebra is just a $\Z$-graded vector
 space
 \[
    V=\bigoplus_{i\in\Z} V_{i}
 \]
 equipped with a bracket operation
 \[
    [\; ,\,]: V_{i} \times V_{j} \longrightarrow V_{i+j}
 \]
 which is super anti-commutative in the sense that
 \begin{equation}
   [x,y]+(-1)^{ij}[y,x]=0, \;\forall \, x\in V_{i}, y\in V_{j}
 \end{equation}
 and satisfies the super Jacobi identity
 \begin{equation} \label{Eq:Jacobi}
   [x,[y,z]] = [[x,y],z]+(-1)^{ij}[y,[x,z]], \; \forall \, x\in V_{i},
   y\in V_{j}.
 \end{equation}
 It is easy to see that the bracket relations defined above do give
 a super Lie algebra structure on $\hat{\g}$. See \cite{DM99} for
 more details on super vector spaces and super Lie algebras.

 \begin{definition}
  A \textit{differential space} is a super vector space $E$ with a
  differential, i.e. an odd endomorphism $d \in End(E)$ satisfying
  $d \circ d = 0$. Endomorphisms of $E$ that commute with the differential $d$
  will be called chain maps or differential space homomorphisms. We use
  $H^{\ast}(E,d)$ to denote the cohomology of $E$ with respect to $d$.
  $(E,d)$ is called \textit{acyclic} if
   \[
     H^{k}(E,d)=
     \begin{cases}
      \R \quad &\text{$k=0$} \\
      0  \quad &\text{$k\neq 0$}
     \end{cases}
   \]
 \end{definition}

 \begin{definition}
  A \textit{homotopy operator} between two chain maps
  $\varphi_{1},\varphi_{2}:E \rightarrow E'$ is an odd linear map
  $h$ (if $E$ and $E'$ are $\Z$ graded, we require h to be of degree
  -1) such that $dh+hd=\varphi_{1}-\varphi_{2}$. if \:
  $\varphi_{1},\varphi_{2}$ are chain homotopic, they induce the
  same map in cohomology.
 \end{definition}

 If we have a Lie group $G$ acting on a differential space $V$, we
 introduce a notion of $G$-differential space.

 \begin{definition}
  A \textit{G-differential space} is a super vector space $V$, together
  with a super Lie algebra homomorphism $\rho: \widehat{\g} \rightarrow
  End(V)$. The horizontal subspace $V_{hor}$ is the space fixed by
  the action of $\g_{-1}$ under $\rho$, the invariant subspace $V^{G}$
  is the space fixed by $\g_{0}$, and the space $V_{basic}$ of
  basic elements is $V_{hor} \cap V^{G}$. It is easy to see that
  $d: V_{basic} \longrightarrow V_{basic}$, i.e. $V_{basic}$ is a
  differential subspace of $V$, we call the the cohomology of
  $(V_{basic},d)$ the \textit{basic cohomology} of the $G$-differential space
  $V$, denoted by $H^{\ast}_{bas}(V)$.

  A \textit{G-differential
  algebra} is a super algebra B with a structure of a
  G-differential space such that $\rho$ takes values in the
  derivation space $Der(B)$ of the algebra.
 \end{definition}

 \begin{example}
   $\widehat{\g}$ is a $G$-differential space with respect to the
   adjoint action of itself.
 \end{example}

 \begin{example}
  Suppose a Lie group $G$ acts on a smooth manifold $M$, i.e. we
  have a group homomorphism $\rho : G \longrightarrow Diff(M)$.
  Then the infinitesimal action $d\rho : \g \longrightarrow
  Vect(M)$ give a representation of the Lie algebra $\g$ of $G$.
  For $\forall \xi \in \g$, let $X_{\xi}$ be the vector field on M
  corresponds $\xi$ under $d\rho$. Let $L_{\xi}$ and $\iota_{\xi}$
  be the Lie derivative and interior product of $X_{\xi}$ in the
  algebra $\Omega^{\ast}(M)$ of smooth differential forms. Then
  $\Omega^{\ast}(M)$ is a $G$-differential algebra with the action
  of $\widehat{\g}= \g_{-1} \oplus \g_{0} \oplus \g_{1}$ by $\iota_{\xi},\ d$
  and $L_{\xi}$.  This example is actually the geometric origin of
  the concept of $G$-differential algebras.
 \end{example}

 \begin{example}
   Choosing a basis $\theta^{1}, \ldots, \theta^{n}$ of $\co{\g}$, we can
   make the exterior algebra $\ext{\co{\g}}$ a $G$-differential algebra by defining:
   \begin{align}
     \iota_{a}\theta^{b}  &= \delta^{b}_{a}, \\
     L_{a} \theta^{b}     &= -f^{b}_{ak}\theta^{k}, \\
     d\theta^{a}          &= -\half f^{a}_{bc}\theta^{b}\theta^{c}
   \end{align}
 \end{example}

  Please note here, I omit the $\wedge$ symbol between two odd
  elements. In my thesis, I will always omit it as long as there
  is no confusion.
  For a compact Lie group $G$, its De Rham cohomology
  $H^{\ast}_{DR}(G)$ coincides with $H^{\ast}(\ext{\co{\g}},d)$
  defined above. So in general $(\ext{\co{\g}},d)$ is not acyclic.

  Similarly, we can define the notion of homomorphism and homotopy in the
  category of $G$-differential spaces.

 \begin{definition}
  A \textit{$G$-homomorphism} between $G$-differential spaces
  $(V_{1},\rho_{1})$ and $(V_{2},\rho_{2})$ is a homomorphism of
  super vector spaces $\phi : V_{1} \rightarrow V_{2}$
  that commutes with the actions of $\widehat{\g}$ on $V_{1}$ and $V_{2}$, i.e.
  for $\forall\, x \in \widehat{\g}, v \in V_{1}, \:
  \phi(\rho_{1}(x)\cdot v)=\rho_{2}(x)\cdot \phi(v)$.
 \end{definition}

 \begin{definition}
  Two $G$-homomorphism $\phi_{1}$, $\phi_{2}$ between two $G$-differential
  spaces $V_{1}$ and $V_{2}$ is called $G$-chain homotopic if there is an odd
  linear map $h: V_{1} \rightarrow V_{2}$ such that $\phi_{1} -
  \phi_{2}= dh+hd,\:  \iota_{a}h+h\iota_{a}=0,\:$ and $
  L_{a}h-hL_{a}=0$.
 \end{definition}

 %%%%%%%%%%%%%%%%%%%%%%%%%%%%%%%%%%%%%%%%%%%%%%%%%%%%%%%%%%%%%%%%%
 %%%%%%%%%% Koszul complex  %%%%%%%%%%%%%%%%%%%%%%%%%%%%%%%%%%%%%%
 %%%%%%%%%%%%%%%%%%%%%%%%%%%%%%%%%%%%%%%%%%%%%%%%%%%%%%%%%%%%%%%%%

 \subsection{Koszul complex} Let $V$ be an
 $n$-dimensional vector space, and let $\wedge(V)$ be the exterior
 algebra of $V$. Koszul algebra $K_{V}$ is the tensor product
 $\wedge(V) \otimes \sym{V}$. The elements $x \otimes 1 \in
 \wedge^1(V) \otimes S^0(V)$ and $1 \otimes x \in \wedge^0(V)
 \otimes S^1(V) $ generate $E_{V}$. The Koszul operator $d_{K}$ is
 defined as the derivation extending the operator on the generators
 given by
  \[
    d_{K}(x \otimes 1)= 1 \otimes x
  \]
  \[
   d_{K}(1 \otimes x)= 0
 \]

 Clearly $d_{K}^2=0$ on generators, and hence everywhere, since
 $d_{K}^2$ is a derivation. In addition, $K_{V}$ can be naturally
 graded by the sum of the natural gradings of $\wedge(V)$ and
 $S(V)$. It is easy to see that $(K_{V},d_{K})$ is an acyclic
 space. Let $x_1, \cdots , x_n$ be a basis of $V$ and define
 \[
   \theta_i:=x_i \otimes 1
 \]
 \[
   z_i := 1 \otimes x_i
 \]

  Then $d=d_{K}$ is expressed in terms of these generators as
 \[
   d\theta_i=z_i
 \]
 \[
    dz_i = 0
 \]

%%%%%%%%%%%%%%%%%%%%%%%%%%%%%%%%%%%%%%%%%%%%%%%%%%%%%%%%%%%%
%%%%%%%%  Weil algebra and G-differentail structure %%%%%%%%
%%%%%%%%%%%%%%%%%%%%%%%%%%%%%%%%%%%%%%%%%%%%%%%%%%%%%%%%%%%%

   \subsection{Weil algebra and its G-differential structures}
   For any Lie algebra $\g$, we define the \textit{Weil algebra}
   $\Weil{\g} = \sym{\co{\g}} \otimes \ext{\co{g}}$. It has a
   natural $G$-differential algebra structure.
   Let $v^{a}=e^{a} \otimes 1$ and $\theta^{a}=1 \otimes e^{a}$.
   They generates the $\sym{\co{\g}}$ and $\ext{\co{\g}}$ in $\Weil{\g}$,
   respectively. We define the action of $\widehat{\g}$ on $\Weil{\g}$
   by
  \begin{align}
      L_{a}v^{b}        &= -f^{b}_{ak} v^{k},   \label{B:begin} \\
      L_{a}\theta^{b}     &= -f^{b}_{ak} \theta^{k},  \\
      \iota_{a}v^{b}    &= 0,         \\
      \iota_{a}\theta^{b} &= \delta^{b}_{a}, \\
      dv^{a}            &= -f^{a}_{bc} \theta^{b} v^{c}, \\
      d\theta^{a}         &= v^{a} - \half
                               f^{a}_{bc}\theta^{b}\theta^{c}
                               \label{B:end}
  \end{align}

   It is easy to see that the horizontal subspace $(\Weil{\g})_{hor} \cong
   \sym{\co{\g}}$ and the basic subalgebra is just the algebra of
   invariant polynomials $\sym{\co{\g}}^{\g}$ on $\g$. In addition,
   since $dv^a=(\theta^{b}L_{b})v^a$ and $d$ is a derivation, we
   conclude that $d\omega=\theta^{b}L_{b}\omega$ for $\forall\, \omega
   \in (\Weil{\g})_{hor}$. So the differential $d$
   actually vanishes on $(W_{\g})_{basic}$. Therefor $H^{\ast}_{bas}(\Weil{\g})
   = (\Weil{\g})_{basic}$.

   We will see below that $\Weil{\g}$ equals the Koszul algebra of
   $\co{g}$ as differential spaces. Hence it is acyclic.

 \begin{proposition} \label{B:Weil-Is-Acyclic}
   $(\Weil{\g},d)$ is an acyclic differential space.
 \end{proposition}

 \begin{proof}
  The easiest way to see this is considering a variable
   change in $\Weil{\g}$. Let
 \begin{equation}
   \widehat{v}^{a}=v^{a} - \half f^{a}_{jk} \theta^{j}\theta^{k}
 \end{equation}
   Extend this naturally to all elements in $\Weil{\g}$.
   Observe that $\widehat{v}^{1}, \ldots , \widehat{v}^{n}, \theta^{1}, \ldots,
   \theta^{n}$ also generates $\Weil{\g}$. The induced
   $G$-differential structure under this set of generators is
 \begin{align}
      L_{a}\widehat{v}^{b}      &= -f^{b}_{ak} \widehat{v}^{k},   \label{A:eqn1} \\
      L_{a}\theta^{b}           &= -f^{b}_{ak} \theta^{k},  \\
      \iota_{a}\widehat{v}^{b}  &= -f^{b}_{ak} \theta^{k},         \\
      \iota_{a}\theta^{b}       &= \delta^{b}_{a}, \\
      d\widehat{v}^{a}                    &= 0,  \label{A:eqn5}\\
      d\theta^{a}               &= \widehat{v}^{a} \label{A:eqn6}
 \end{align}

   Equation~\eqref{A:eqn5} and~\eqref{A:eqn6} mean that $\Weil{\g}$
   is a Koszul complex generated by $\widehat{v}^{a}$ and $\theta^{a}$.
   So it is acyclic.

 \end{proof}

   We will find that this different presentation of the $G$-differential structure
   on $W(\g)$ in~\eqref{A:eqn1}---~\eqref{A:eqn6} is useful in next chapter.
   So we have the following definition.
  \begin{definition}
    Let $W^{K}_{\g}= S(\co{\g}) \otimes \ext{\co{\g}}$ with
    generators $\widehat{v}^{a}$ and $\theta^{a}$, and define
    its G-differential structure by~\eqref{A:eqn1}---~\eqref{A:eqn6}.
    We call it Koszul G-differential structure of $\Weil{\g}$.
  \end{definition}

  By definition, $\Weil{\g}$ and $W^{K}_{\g}$ are isomorphic as
  G-differential algebras. The isomorphism is just the super
  variable change.
  \begin{equation} \label{D:varchange}
   \begin{array}{rccl}
    \tau_{0}: & W^{K}_{\g} & \longrightarrow &  \Weil{\g} \\
    \     &   \widehat{v}^{a} & \mapsto & v^{a}- \half f^{a}_{jk}\theta^{j}\theta^{k} \\
    \     & \theta^{b}& \mapsto & \theta^{b}
   \end{array}
   \end{equation}

   %%%%%%%%%%%%%%%%%%%%%%%%%%%%%%%%%%%%%%%%%%%%%%%%%%%%%%
   %%%%%% G-differential structure of Cl(g) %%%%%%%%%%%%%
   %%%%%%%%%%%%%%%%%%%%%%%%%%%%%%%%%%%%%%%%%%%%%%%%%%%%%%

\subsection{G-differential structure of Clifford algebra}
    Suppose a Lie algebra $\g$ is equipped with an $ad(\g)$-invariant inner
    product $B(\, ,\,)$, we call this type of Lie algebra
    \textit{quadratic Lie algebra}.
    Semisimple Lie algebras and compact Lie algebras are all
    quadratic. However, not every Lie algebra allows an invariant
    inner product.
    From now on, we always assume $\g$ is quadratic.
    Let $\Cl{\g}$ be the \textit{Clifford algebra}
    of $\g$, i.e. the quotient of the tensor algebra $T(\g)$ by the ideal
    generated by all $2x \otimes x - B(x,x)$ with $x \in \g$. It
    inherits a natural $\Z_{2}$-grading from the tensor algebra
    and a filtration,
    \[
      \R=Cl^{(0)}(\g) \subset Cl^{(1)}(\g) \subset \ldots
    \]
    So $\Cl{\g}$ is a filtered super algebra. The associated
    graded algebra $Gr^{\ast}(\Cl{\g})$ is isomorphic to
    $\ext{\g}$. Let elements in $\g$ have odd grading and
    define a bracket operation in $\g^{odd}:=\g \oplus \R\mathfrak{c}$ by
    \begin{align}
      [x,y]_{odd}  &=B(x,y)\mathfrak{c},\: \forall\ x,y \in \g  \notag  \\
      [x,\mathfrak{c}]_{odd} &=0, \: \forall\ x \in \g \notag
    \end{align}

    Here $\mathfrak{c}$ is an even element. We use a subscript $odd$
    to distinguish this bracket from the original Lie bracket of $\g$.
    Then we can look on $\Cl{\g}$ as the universal enveloping
    algebra of $\g\oplus \F\mathfrak{c}$ modulo $\mathfrak{c}=1$.
    The odd bracket of $\g^{odd}$ can be extended naturally to
    $\Cl{\g}$. The general formula is :
    \begin{align} \label{odd_bracket}
       & \quad\: [x_{a_1} \ldots x_{a_k}, x_{b_1}, \ldots, x_{b_l}]_{odd}
         \notag \\
      &= \underset{1\leq j\leq l}{\sum_{1\leq i \leq k }}
             (-1)^{k-i-j-1}[x_{a_i},x_{b_j}]_{odd} \:
             x_{a_1} \ldots \widehat{x}_{a_i} \ldots x_{a_k}
               x_{b_1} \ldots \widehat{x}_{b_j} \ldots x_{b_l}
   \end{align}

    The Poincar\'{e}-Birkhoff-Witt map in this case is just the
    anti-symmetrization, which is also called Chevalley \textit{quantization map}.
   \[
      q: \ext{\g} \longrightarrow \Cl{\g}
   \]
   \[
      q(x_{1}\wedge \ldots \wedge x_{k})=
         \sum_{\sigma \in S_{k}}
               (-1)^{sgn(\sigma)}
                    x_{\sigma(1)} \ldots x_{\sigma(k)}
   \]
    Next, we make $\Cl{\g}$ a $G$-differential algebra by first defining
    the action of $\widehat{\g}$ on a basis of $\g$ and then extend it by derivations
    to the whole $\Cl{\g}$ as follows.
    Suppose $e_{1}, \ldots, e_{n}$ is an orthonormal basis of $\g$ with respect to
    the invariant metric $B(\, ,\,)$. If $[e_{a},e_{b}]=f^{c}_{ab}e_{c}$,
    then from the fact
    $B([x_{a},x_{b}],x_{c})+B(x_{b},[x_{a},x_{c}])=0$, we can
    easily see that $f^{c}_{ab}=-f^{b}_{ac}$. So under this basis, the
    structure constant $f_{ab}^{c}$ is totally anti-symmetric with respect to
    its indices. For convenience, we use $f_{abc}$ to denote $f^{c}_{ab}$
    under an orthonormal basis.
    In addition, we can naturally identify $\co{\g}$ and $\g$ using
    $B(\, ,\,)$. Now, we let $\widehat{\g}$ act on $\Cl{\g}$ by
   \begin{align}
     \iota_{a}(e_{b}) &= [e_{a},e_{b}]_{odd}= \delta_{ab} \\
     L_{a}(e_{b})     &= [e_{a},e_{b}]=[-\half f_{aij}e_{i}e_{j},e_{b}]_{odd}
                           =f_{abc}e_{c} \\
     d(e_{a})         &=
                         [-\frac{1}{6}f_{ijk}e_{i}e_{j}e_{k},e_{a}]_{odd}
                       = -\half f_{aij}e_{i}e_{j}
   \end{align}

   \begin{proposition}
   The above definition gives a $G$-differential algebra structure on $\Cl{\g}$ and
   $(\Cl{\g},d^{Cl})$ has trivial cohomology in all filtration degrees
   for any quadratic Lie algebra $\g$.
   \end{proposition}

   \begin{proof}
    See \cite{AMM98} chapter3.
   \end{proof}

   So generally speaking, the exterior algebra $\ext{\co{\g}}$ and the
   Clifford algebra $\Cl{\g}$ are not isomorphic as differential
   spaces (unless $\g$ is abelian).

 %%%%%%%%%%%%%%%%%%%%%%%%%%%%%%%%%%%%%%%%%%%%%%%%
 %%%%%%%% Noncommutative Weil algebra %%%%%%%%%%%
 %%%%%%%%%%%%%%%%%%%%%%%%%%%%%%%%%%%%%%%%%%%%%%%%

 \subsection{Noncommutative Weil algebra}
    In \cite{AMM98}, A.Alekseev and E.Meinrenken defines a
    noncommutative version of Weil algebra, they call it \textit{
    noncommutative Weil algebra}.
  \begin{definition}
    For a quadratic Lie algebra $\g$ with invariant metric $B(\; , \:)$,
    if $\widehat{u}_1, \ldots , \widehat{u}_n$ is a basis of $\g$ and
    $\xi_{1}, \ldots ,\xi_{n}$ be the corresponding basis of another copy $\g_{1}$ of $\g$,
    the noncommutative Weil algebra $\NWeil{\g}$ of $\g$ is defined to the
    the quotient of the tensor algebra $T(\g \oplus \g_{1})$ by
    the following relations:
    \begin{align}
      \widehat{u}_a \otimes \widehat{u}_b -\widehat{u}_b \otimes \widehat{u}_a
                                 &= [\widehat{u}_a,\widehat{u}_b]_{\g} , \notag \\
      \xi_{a} \otimes \xi_{b} - \xi_{b} \otimes \xi_{a} &= B(\xi_{a},\xi_{b}) , \notag \\
      \widehat{u}_a \otimes \xi_{b} - \xi_{b} \otimes \widehat{u}_a   &= [\xi_{a},
      \xi_{b}]_{\g_{1}} \notag
    \end{align}
  \end{definition}
 %%%%%%%% define the Koszul structure %%%%%%%
  We can define the $G$-differential algebra structure on the noncommutative Weil
  algebra of $\g$ as follows:
  Suppose $\widehat{u}_{a}$ and $\xi_{a}$ have degree $2$ and $1$, respectively.
  Then define:
  \begin{align}
      L_{a}\widehat{u}_{b}         &= f_{abc} \widehat{u}_{c},    \\
      L_{a}\xi_{b}                 &= f_{abc} \xi_{c},  \\
      \iota_{a}\widehat{u}_{b}     &= f_{abc} \xi_{c},         \\
      \iota_{a}\xi_{b}             &= \delta_{ab}, \\
      d\widehat{u}_{a}             &= 0, \\
      d\xi_{a}                     &= \widehat{u}_{a}
  \end{align}

  If we identify $\co{\g}$ and $\g$ using the invariant metric on
  $\g$, it is easy to see that the $G$-differential structure of
  $\NWeil{\g}$  defined above corresponds exactly to the $G$-
  differential structure of $W^{K}_{\g}$.
  So we use $\NWeil{\g}^{K}$ to denote $\NWeil{\g}$ with this
  $G$-differential structure.

  \begin{remark}
   The naive extension of the map $v^a \longrightarrow u_a$,
  $\theta_a \longrightarrow \xi_a$ by symmetrization is not
  a $G$-differential algebra homomorphism from $W^{K}_{\g}$ to
  $\NWeil{\g}^{K}$. It is only a linear map.
  \end{remark}

  In addition, similar to Weil algebra, we can define a super variable change $\tau_{1}$:
  \begin{align}
    \widehat{u}_{a}   & = u_{a} - \half f_{abc} \xi_{b}\xi_{c}  \notag \\
    \xi_{a} & = \xi_{a} \notag
  \end{align}
   Notice $u_{a}$ and $\xi_{a}$ also generate
   $\NWeil{\g}$, and they commute. In fact
  \begin{align}
     [u_a,\xi_b]&=[\widehat{u}_a + \half f_{aij}\xi_i\xi_j,\xi_b]
                 = f_{abi}\xi_{i} + \half f_{aij}\xi_i \delta_{jb} -
                  \half f_{aij}\xi_j \delta_{ib} \\ \notag
          \:    &= f_{abi}\xi_{i} + \half f_{aib}\xi_i - \half f_{abj}\xi_j
                 = 0  \notag
  \end{align}

   Under this variable change, the induced $G$-differential structure on
  $\NWeil{\g}$ defined on generators $u_{a}$ and $\xi_{a}$ is:
   \begin{align}
      L_{a}u_{b}         &= f_{abc} u_{c},    \\
      L_{a}\xi_{b}       &= f_{abc} \xi_{c},  \\
      \iota_{a}u_{b}     &= 0,         \\
      \iota_{a}\xi_{b}   &= \delta_{ab}, \\
      du_{a}             &= -f_{abc} \xi_{b}u_{c}, \\
      d\xi_{a}           &= u_{a} - \half f_{abc}\xi_{b}\xi_{c}
   \end{align}

    Obviously, this form of $G$-differential structure on $\NWeil{\g}$
    corresponds to the $G$-differential structure of $\Weil{\g}$ defined
    in~\eqref{B:begin}--~\eqref{B:end}. When we write $\NWeil{\g}$ without
    the capital $K$, we mean its generators are $u_a$ and $\xi_a$.

  \begin{remark}
    It is not hard to see that the bracket $[\widehat{u}_{a},\widehat{u}_{b}]$
    corresponds exactly to $[u_a,u_b]$ under the super variable change.
    Thus we can think of $u_{a}$ being a basis of a Lie algebra $\g_{0}$
    which isomorphic to $\g$. Then as vector spaces, $\NWeil{\g} = \NWeil{\g}^{K} = U(\g_{0})
    \otimes Cl(\g_{1})$.
  \end{remark}

\subsection{Notation summary}

  To avoid confusions, I summarize different notations of basis elements of $\g$ or $\co{\g}$
  we use in this paper so far.
  \begin{equation}
    \begin{array}{cccccc}
     \g  & \co{\g} & \sym{\co{\g}} & \ext{\co{\g}} & \U{\g} & \Cl{\g} \\
     e_{a} & e^{a} & v^{a},\widehat{v}^a & \theta^{a} & u_{a}, \widehat{u}_a & \xi_{a}
    \end{array}
  \end{equation}

%%%%%%%%%%%%%%%%%%%%%%%%%%%%%%%%%%%%%%%%%%%%%%%%%%%%%%%%%%%%%%%%%%%%%
 %%%%%%%%%  Duflo Isomorphism and Quantization map %%%%%%%%%%%%%%%%%%%
 %%%%%%%%%%%%%%%%%%%%%%%%%%%%%%%%%%%%%%%%%%%%%%%%%%%%%%%%%%%%%%%%%%%%%

  \subsection{Duflo isomorphism and Quantization map} Every Lie algebra $\g$ has two
   associated algebras: the symmetric algebra $\sym{\g}$, generated
   by $\g$ with relations $xy-yx=0$, and the universal enveloping
   algebra $\U{g}$, generated by $\g$ with relations $xy-yx=[x,y]$.
   There is a natural map between the two algebras,
   \[
     \chi:\sym{\g} \rightarrow \U{g}
   \]
   given by taking a monomial $x_1 x_2 \ldots x_n$ in $S(\g)$ and
   averaging over the product in $\U{g}$ of all the $x_{i}$ in all
   possible orders. By the Poincar\'{e}-Birkhoff-Witt(PBW) theorem,
   $\chi$ is an isomorphism of vector spaces and $\g$-modules. $\chi$
   is not an algebra isomorphism, even restricted to the invariant
   subspaces of both sides.  in \cite{Duf77}, M. Duflo modifies
   $\chi$ a little and gives an algebra isomorphism
    \[
      \Upsilon: \sym{\g}^{\g}  \rightarrow  \U{g}^{\g}
    \]
   where
  \begin{align}
      \Upsilon &=\chi \circ \partial_{j^{\half}} \\
      j^{\half}(x) &= det^{\half} \left( \frac{\sinh(\half
      ad_{x})}{\half ad_{x}} \right)
             = det^{\half} \left( \frac{e^{ad(x)/2}-e^{-ad(x)/2}}{ad(x)}
                           \right)
  \end{align}

  The notation $ad_{x}$ is the adjoint action of an element $x$ on
  $\g$. The $\partial_{j^{\half}}$ means to consider $j^{\half}(x)$
  as a power series on $\g$ and so we can think of it as an
  infinite-order differential operator on $\co{g}$, which we can
  then apply to polynomials on $\co{g}$ ($\equiv$ elements in
  $\sym{\g}$). $j^{\half}$ is an important function in the theory of
  Lie algebras. Its square, $j(x)$, is the Jacobian of the
  exponential mapping from $\g$ to its Lie group $G$ when $\g$ is unimodular.
  The map $\Upsilon$ here is called Duflo isomorphism. When Lie algebra $\g$
  is semisimple, Duflo isomorphism coincides with the Harish-Chandra
  isomorphism. The general proof in \cite{Duf77} is highly
  nontrivial, it used certain facts about finite-dimensional Lie
  algebras which follow only from the classification theory. In
  addition, Duflo isomorphism is also true in the case of Lie super
  algebras, which is proved by M. Kontsevich in \cite{Ko97}.

  In addition, we have another way of understanding the Duflo map
  which is more convenient to use in many cases. We can identify
  $\sym{\g}$ with the space of distributions which are supported at the origin of $\g$.
  In fact, on any vector space V with a
  basis $\{ e_a \}$, let $x_a$ be the coordinate function of $e_a$, the algebra
  $\mathcal{E}_{\textbf{0}}(V)$ of distributions with support at the origin on V is
  canonically isomorphic to its symmetric algebra $S(V)$, by identifying each
  basis element $e_a$ with $-\pardev{}{x_a} \delta_{0}$,  Where $\delta_{0}$ is
  the \textit{Dirac delta function} at the
  origin of $V$.
         \[
            e_{a} = -\pardev{}{x_{a}} \delta_{0},
         \]
 The algebra structure on $\sym{V}$ corresponds to the
 convolution $\ast_{V}$ of distributions (see \cite{Fried98} for details).
 In addition, any smooth function $f(x)$ can act on
 $\mathcal{E}_{\textbf{0}}(V)$ by
  \[
       \langle f(x) \cdot D, \phi(x) \rangle = \langle D, f(x)
       \phi(x)
       \rangle ,
       \;\; \forall \: D \in \mathcal{E}_{\textbf{0}}(V), \phi(x) \in C^{\infty}_0(V).
  \]

  Similarly for $\U{\g}$, we can identify
  the generator $u_a$ with a distribution on $G$ which is supported
  at the identity element $\textbf{1}$ of $G$. We have
  \[
    u_a = \left. \frac{d}{dt} \right|_{t=0}
                 \delta_{\exp(te_{a})}
  \]
  In this way, $\U{\g}$ can be identified algebraically with the space
  of distributions $\mathcal{E}_{\textbf{1}}(G)$ on $G$ supported at the
  identity $\textbf{1}$ of $G$. The product in $\U{\g}$
  corresponds to the convolution product $\ast_{G}$ of distributions in
  $\mathcal{E}_{\textbf{1}}(G)$.

  If we view $\sym{\g}$, $\U{\g}$ and the action of $j^{\half}(x)$ as
  above described way,
  the PBW map $\chi$ is interpreted as pushing forward
  distributions on $\g$ to distributions on $G$ by the exponential
  map of $\g$. In fact, this is the point of view adopted by M. Duflo
  in his original work \cite{Duf77}. So we can write Duflo's theorem as follows:
  \begin{equation}
      \Upsilon(\eta)=\exp_{\ast}(j^{\half}(x)\cdot \eta), \:
         \forall \ \eta \in \mathcal{E}_{\textbf{0}}(\g)=\sym{\g}
  \end{equation}
  \begin{equation}
      \Upsilon(\eta_1 \ast_{\g} \eta_2)= \Upsilon(\eta_1) \ast_{G}
       \Upsilon(\eta_2), \: \forall \ \eta_1,\eta_2 \in
       \sym{\g}^{\g}
  \end{equation}

   In fact, these two ways of understanding Duflo map are Fourier transform
   of each other in a canonical way.
  
  More recently, A.Alekseev and E.Meinrenken establishes in
  \cite{AMM98}
  an interesting isomorphism $\mathcal{Q}$ between $\ext{\co{g}}
  \otimes \sym{\co{g}}$ and $\Cl{g} \otimes \U{g}$ for any quadratic
  Lie algebra $\g$. They called $\mathcal{Q}$ the
  \textit{quantization map}. If we identity $\co{\g}$ and $\g$ using
  the metric, when restricted to the symmetric algebra $1 \otimes
  \sym{\co{g}}$, the quantization map $\mathcal{Q}$ becomes the
  usual Duflo isomorphism from $\sym{\g}$ to $\U{g}$. When
  restricted to the exterior algebra $\wedge(\co{g}) \otimes 1$,
  $\mathcal{Q}$ becomes the antisymmetrization map $q$ from
  $\wedge(\g)$ to $\Cl{\g}$. However, $\mathcal{Q}$ is not just the
  direct product $Duf \otimes q$. It has the more complicated form
   \[
      \mathcal{Q}=(Duf \otimes q) \circ exp \left(
                 \half T_{ab}(x) \iota_{a}\iota_{b} \right)
              =  (\chi \otimes q) \circ \partial_{j^{\half}(x)} \circ
                  exp \left( \half T_{ab}(x)\iota_{a}\iota_{b} \right)
   \]

  where  $T_{ab}$ is a certain anti-symmetric tensor field on
  $\g$, i.e. an $\wedge^{2}(\g)$ value function on $\g$. It
  was obtained by Etingof-Varchenko \cite{EtVar1998} as a solution
  of the classical dynamical Yang-Baxter equation.
  \[
     T_{ab}(x):= (\ln(j)'(ad_{x}))_{ab}
  \]
  When acting on the elements in Weil algebra $W(\g)$, $T_{ab}(x)$
  is treated as a differential operator like $j^{\half}(x)$,
  $\iota_a\iota_b$ is understood as contractions of odd variables.

  We will see in Chapter 4 that, the quantization map $\mathcal{Q}$
  can be understood as the super Duflo map of a super Lie algebra.
  \newline

  %%%%%%%%%%%%%%%%%%%%%%%%%%%%%%%%%%%%%%%%%%%%%%%%%%%%%%%%% end of chapter
 %Definitions and Preliminary Facts

\large
\section{\textbf{Supergeometric Interpretation of the Quantization Map}}
\normalsize
  %Introduce the Lie algebra of the super Lie group
  Suppose $G$ is a Lie group with Lie algebra $\g$. Set $TG[1]=G \times
  \g$ which we look upon as a Lie group with multiplication given
  as follows
  \[
   (g,X)(h,Y)=(gh,Ad_{h^{-1}}(X)+Y), \quad  \forall \, g,h \in G, X,Y \in \g
  \]

  $TG[1]$ is just the odd tangent bundle of $G$ with the natural group
  structure induced from the group structure of $G$. The Lie
  algebra, $T\g[1]$, of $TG[1]$ is $\g \times \g$ with bracket given by
  \[
   [(X,Y),(X',Y')]=([X,X'],[Y,X']+[X,Y']), \quad  \forall X,X',Y,Y'
   \in \g
  \]

   In $T\g[1]$, $\g \times 0$ is a Lie subalgebra isomorphic to $\g$
   and $0 \times \g$ is an abelian Lie subalgebra. When $\g$ is
   a quadratic Lie algebra, we can make a central extension of the
   Lie algebra $T\g[1]$ using its invariant metric $B(\:\: ,\, )$,
   \begin{equation} \label{Intro:extension}
     0 \longrightarrow \R \longrightarrow \widetilde{T\g[1]}
     \longrightarrow T\g[1] \longrightarrow 0
   \end{equation}

   And the Lie bracket of $0 \times \g$ in $T\g[1]$ becomes
   \[
       [(0,X),(0,Y)]=B(X,Y) \mathfrak{c}
   \]

   where $\mathfrak{c}$ is a central element in $\widetilde{T\g[1]}$.
   It generates the $\R$ in~\eqref{Intro:extension}.

  %%%%%%%%%%%%%%%%%%%%%%%%%%%%%%%%%%%%%%%%%%%%%%%%%%
  %%%%%%%% Relations between S(Tg[1]) and W(g) %%%%%
  %%%%%%%%%%%%%%%%%%%%%%%%%%%%%%%%%%%%%%%%%%%%%%%%%%

 \subsection{Relations between $S(\widetilde{T\g[1]})$ and $\Weil{\g}$ and $\Weil{\g}^{K}$ }
   If we identify $\co{\g}$ and $\g$ via the invariant metric $B(\:,\,)$,
   $\Weil{\g}^{K}$ is then identified with $S(\widetilde{T\g[1]})/< \mathfrak{c}=1>$.
   In fact, we can think of $W^{K}_{\g}$ as the space of differential forms
   on $\g$ with polynomial coefficients.
   Let $\sigma_{0}: S(\widetilde{T\g[1]}) \longrightarrow \NWeil{\g}$
   be a map which make the following diagram commute:
  \[
    \begin{diagram}
        \node{\text{$W_{\g}^{K}$}}
             \arrow{se,l}{\text{$\tau_{0}$}}
        \node{\text{$S(\widetilde{T\g[1]}$)}}
             \arrow{w,tb}{\text{$\mathfrak{c}=1$}}{B(,)}
             \arrow{s,r}{\text{$\sigma_{0}$}}    \\
       \node[2]{\text{$W_{\g}$}}
    \end{diagram}
  \]

   Since
   $\tau_{0}(\widehat{v^{a}}) = v^{a}- \half
                        f_{abc}\theta^{b}\theta^{c}$,
   if $e_{a}$ and $\overline{e_{a}}$ are basis of
   $\widetilde{T\g[1]}^{even}$ and $\widetilde{T\g[1]}^{odd}$
   respectively, $\sigma_{0}$ must be:
   \begin{align}
     \sigma_{0}(e_{a}) &= v^{a} -
                              \half f_{abc}\theta^{b}\theta^{c} \notag \\
     \sigma_{0}(\overline{e_{a}}) &= \theta^{a} \notag
   \end{align}

    The following lemma tells us how the super Lie algebra
  structure of $\widetilde{T\g[1]}$ is related to the
  $G$-differential structure of $\Weil{\g}$ under the map $\sigma_{0}$.

  \begin{lemma}
    The map $\sigma_{0}: S(\widetilde{T\g[1]}) \longrightarrow
    \Weil{\g}$ satisfies:
    \begin{align}
     \sigma_{0}([e_{a},e_{b}])
               &= L_{a}(\sigma_{0}(e_{b}))  \label{3:equ1}\\
     \sigma_{0}([e_{a},\overline{e_{b}}])
               &= L_{a}(\sigma_{0}(\overline{e_{a}})) \\
     \sigma_{0}([\overline{e_b},e_a])
               & = \iota_{b}(\sigma_{0}(e_{a}))  \\
     \sigma_{0}([\overline{e_{a}},\overline{e_{b}}])
               &= \iota_a(\sigma_{0}(\overline{e_{b}}))
    \end{align}
  \end{lemma}

  \begin{proof}
    Here I only give the proof of~\eqref{3:equ1}, the proof of the other
    two equations are similar. By the definition of $\sigma_{0}$, we have
    \[
     \sigma_{0}([e_{a},e_{b}])
        = \sigma_{0}(f_{abc}e_{c})
        = f_{abc}(v^{c} -
                   \half f_{cpq}\theta^{p}\theta^{q})
    \]
    On the other hand,
    \begin{align}
     L_{a}(\sigma_{0}(e_{b}))
        &= L_{a}(v^{b}-\half
          f_{bij}\theta^{i}\theta^{j})
         = f_{abc}v^{c} -
              L_{a}(\half f_{bij}\theta^{i}\theta^{j}) \notag \\
     \  &= f_{abc}v^{c} -
            \half f_{bij}f_{aik}\theta^{k}\theta^{j}-
            \half f_{bij}f_{ajk}\theta^{i}\theta^{k}   \notag \\
     \  &= f_{abc}v^{c} -
            \half f_{abk}f_{kij}\theta^{i}\theta^{j} \notag
    \end{align}
    The last equality is because of the Jacobi identity. So the~\eqref{3:equ1} holds.
   \end{proof}

  By the above lemma, we have the following correspondence:
 \[
   \begin{array}{cc}
    \sigma_{0}(e_{a}) \longleftrightarrow e_{a}
    & L_{a} \longleftrightarrow [\, e_{a},\;] \\
      \theta_{a} \longleftrightarrow \overline{e_{a}}
    & \iota_{a} \longleftrightarrow [\overline{e_{a}},\:]
   \end{array}
 \]

   Since $\sigma_0$ is an algebra isomorphism, we can easily see the following.
  \begin{proposition}
   The basic complex of $\Weil{\g}^{K}$ corresponds exactly to the
   invariant symmetric tensor of $\widetilde{T\g[1]}$ under
   $\sigma_{0}$. Therefore, $S(\widetilde{T\g[1]})^{inv}\cong S(\g)^{\g}$.
  \[
   \sigma_{0}: S(\widetilde{T\g[1]})^{inv}   \longleftrightarrow (\Weil{\g})_{basic}
  \]
  \end{proposition}

  \  \

  In addition, it is easy to see that
  $\NWeil{\g} = U(\widetilde{T\g[1]})/<\mathfrak{c}=1>$ by definition.
  So we get a similar diagram for the noncommutative
  Weil algebra $\NWeil{\g}$ and $\NWeil{\g}^{K}$.
 \[
  \begin{diagram}
   \node{\text{$U(\widetilde{T\g[1]}$)}}
             \arrow{e,t}{\text{$\mathfrak{c}=1$}}
             \arrow{s,l}{\text{$\sigma_{1}$}}
   \node{\text{$\NWeil{\g}^{K}$}}
             \arrow{sw,r}{\text{$\tau_{1}$}} \\
   \node{\text{$\NWeil{\g}$}}
  \end{diagram}
 \]
  Where $\tau_{1}$ is the super variable change. The same argument
  as the preceding lemma shows that $\sigma_1$ maps $U(T\g[1])^{inv}$
  isomorphically onto $(\NWeil{\g})_{basic}$.

 %%%%%%%%%%%%%%%%%%%%%%%%%%%%%%%%%%%%%%%%%%%%%%%%
 %%%%%%%% Duflo isomorphism of Tg[1] %%%%%%%%%%%%
 %%%%%%%%%%%%%%%%%%%%%%%%%%%%%%%%%%%%%%%%%%%%%%%%

\subsection{Duflo isomorphism of $\widetilde{T\g[1]}$}
  In chapter two , we have introduced the Duflo isomorphism of a
  Lie algebra and it actually makes sense for any super Lie algebras.
  Now, let us examine the Duflo isomorphism of $\widetilde{T\g[1]}$.

  Let $\{e_{1},\ldots,e_{n},\overline{e_{1}},\ldots,\overline{e_{n}},\mathfrak{c}\}$
   be a basis $\widetilde{T\g[1]}$. By definition, $ad(e_{a})$ and $ad(\overline{e_{a}})$ are represented
   by the following matrices :
  \begin{align}
          \       &\quad \,\,\; even \;\;\; odd \;\;\: \mathfrak{c}  \notag \\
     ad(e_{a})&= \left(
                       \begin{matrix}
                          M_a & 0   & 0 \, \\
                          0   & M_a & 0 \, \\
                          0   & 0   & 0 \,
                       \end{matrix}
                    \right)
                       \begin{matrix}
                         even \\
                         odd \\
                         \mathfrak{c}
                       \end{matrix}   \notag \\
    ad(\overline{e_{a}})&= \left(
                           \begin{matrix}
                                0   & \; 0\;\; & 0 \, \\
                               M_a  & \; 0\;\; & 0 \, \\
                       \epsilon_{a} & \; 0\;\; & 0 \,
                           \end{matrix}
                    \right)
                       \begin{matrix}
                         even \\
                         odd \\
                         \mathfrak{c}
                       \end{matrix}   \notag
  \end{align}
  Where $(M_{a})_{ij}=f_{aji}$ and
  $\epsilon_{a}=(0,\ldots,0,\stackrel{a}{1},0,\ldots,0)$.

  Then it is easy to see that for any element $x \in
  \widetilde{T\g[1]}$, $Tr(ad^{k}(x))\equiv 0$ for $\forall k \in \N$
  (here $Tr$ is the super trace, see \cite{DM99} for definition).

  \[
     j^{\half}(x)= \det \left( \frac{\sinh(\half
                     ad_{x})}{\half ad_{x}} \right)
                 = \exp{\left(\sum_{k=1}^{\infty}b_{2k}Tr(ad^{k}(x))\right)}
                 = 1
  \]
    Where $b_{2k}$'s are modified Bernoulli numbers defined by the
    power series expansion
  \begin{equation} \label{eq:Bernoulli}
     \sum_{k=0}^{\infty}b_{2k}x^{2k}=\half
     \ln{\left( \frac{\sinh{\frac{x}{2}}}{\frac{x}{2}} \right)}.
  \end{equation}

  So we have proved the following.

  \begin{proposition}
  the Duflo isomorphism for $\widetilde{T\g[1]}$ is just the
  (super) symmetrization map.
  \end{proposition}

%%%%%%%%%%%%%%%%%%%%%%%%%%%%%%%%%%%
%%%%%% Main Theorem %%%%%%%%%%%%%%%
%%%%%%%%%%%%%%%%%%%%%%%%%%%%%%%%%%%

 \subsection{The Main Theorem}

   I want to show that the quantization map between Weil algebra and
  noncommutative Weil algebra is essentially equivalent to the Duflo map for
  the super Lie algebra $\widetilde{T\g[1]}$. I put this in the main theorem
  as follows.

  \begin{main theorem}
  The quantization map $\mathcal{Q}$
  is the (super) Duflo isomorphism of the super Lie algebra
  $\widetilde{T\g[1]}$, i.e. the following diagram commutes.

 % Insert Diagram 1
 \begin{equation} \label{Main_Diag_1}
  \begin{diagram}
   \node{\text{$W_{\g}^{K}$}}
        \arrow{se,b}{\text{$\tau_{0}$}}
   \node{\text{$S(\widetilde{T\g[1]}$)}}
        \arrow{w,tb}{\text{$\mathfrak{c}=1$}}{B(,)}
        \arrow[2]{e,tb}{\text{Duflo}}{\text{vector space homo.}}
        \arrow{s,r}{\text{$\sigma_{0}$}}
   \node[2]{\text{$U(\widetilde{T\g[1]}$)}}
             \arrow{e,t}{\text{$\mathfrak{c}=1$}}
             \arrow{s,l}{\text{$\sigma_{1}$}}
   \node{\text{$\NWeil{\g}^{K}$}}
             \arrow{sw,b}{\text{$\tau_{1}$}}  \\
   \node[2]{\text{$W_{\g}$}}
             \arrow[2]{e,tb}{\text{$\mathcal{Q}$}}{\text{vector
             space homo.}}
   \node[2]{\text{$\NWeil{\g}$}}
 \end{diagram}
 \end{equation}

 \end{main theorem}

  \ \

 \begin{remark} \label{commute}
  When we apply $\mathcal{Q}$ to elements in $W^{K}_{\g}$, we use
  the metric $B( \;\: , \;)$ to identify $W_{\g}$ with
  $\sym{\g} \otimes \ext{\g}$.
 \end{remark}

   \  \

 \begin{corollary}
  The basic cohomology of $\Weil{\g}$ is algebraically isomorphic to the
  basic cohomology of $\NWeil{\g}$.
 \end{corollary}

 \begin{proof}
   Since Duflo map is an algebra isomorphism when restricted to the
   invariants and $\sigma_{0},\sigma_{1}$ in the diagram are also algebra
   homomorphisms, so the quantization map $\mathcal{Q}|_{(W_{\g})_{basic}}$
   is an algebra homomorphism by the commutativity of the diagram above.
   So we have the following diagram.
  % Insert Diagram 2
 \begin{equation} \label{Main_Diag_2}
  \begin{diagram}
   \node{\text{$\sym{\widetilde{T\g[1]}}^{inv}$}}
          \arrow[3]{e,tb}{\text{Duflo}}{\text{algebra homo.}}
          \arrow{s,l}{\text{$\sigma_{0}$}}
   \node[3]{\text{$U(\widetilde{T\g[1]})^{inv}$}}
          \arrow{s,r}{\text{$\sigma_{1}$}}     \\
   \node{\text{$(W_{\g})_{basic}$}}
          \arrow[3]{e,tb}{\text{$\mathcal{Q}$}}{\text{algebra homo.}}
   \node[3]{\text{$(\NWeil{\g})_{basic}$}}
  \end{diagram}
 \end{equation}

  Notice the differential $d$ is actually trivial on
  $(W_{\g})_{bas}$, so we have $H^{\ast}_{bas}(W_{\g})= (W_{\g})_{basic}= S(\co{\g})^{\g}$.
  Therefore the quantization map $\mathcal{Q}$ induces an algebra homomorphism in
  the basic cohomology.
 \end{proof}

 \begin{remark}
\end{remark}
\begin{enumerate}
    \item The action of the quantization map $\mathcal{Q}$ on $H^{\ast}_{bas}(W_{\g})$
  is just the usual Duflo map for the Lie algebra $\g$. Although
  this is a little disappointing, we will see in the last chapter that
  the quantization map $\mathcal{Q}$ itself has an interesting diagrammatic interpretation.

    \item  The above corollary can also be derived from a very strong theorem
  proved in \cite{AM03} which asserts that any $G$-differential space
  homomorphism from the Weil algebra of $\g$ to a locally free
  $G$-differential algebra (may not be commutative) always induces an algebra homomorphism
  in basic cohomology. \newline
 \end{enumerate}

  %%%%%%%%%%%%%%%%%%%end of Chapter 3
 %Supergeometric interpretations of Quantization Map

\large
\section{\textbf{Proof of the Main Theorem}}
\normalsize

   In this chapter, I will present an algebraic proof of the the
   main theorem in this thesis. To do that, we need to investigate the
   structure of the spin representation of the spin group. The techniques
   in the proof are described by A.Alekseev and E.Meinkenren in
   \cite{AM03}\cite{AMM98}\cite{AM02}.
   So I will quote some theorems directly from their papers without
   giving the proof since they are quite complicated themselves.
   First, let's look at the spin representation.

%%%%%%%%%%%%%%%%%%%%%%%%%%%%%%%%%%%%%%%%%%%%%%%%%%%%%%%%%%%%%%%%
%%%%%%% Spin Group of Cl(g) and spin representation  %%%%%%%%%%%%%%%
%%%%%%%%%%%%%%%%%%%%%%%%%%%%%%%%%%%%%%%%%%%%%%%%%%%%%%%%%%%%%%%%

\subsection{Spin group of $\g$ and spin representation}
  In Chapter two, I have introduced the Clifford algebra of a
  quadratic Lie algebra $\g$ and defined a $G$-differential structure
  on it. In this section, I will discuss an important object
  $Spin(\g)$ in $\Cl{\g}$ which will be used in the proof of the
  main theorem. People
  can find more details of Clifford algebra and $Spin(\g)$ and its
  representations
  in \cite{Lawson89}, \cite{GoodNolan98},
  \cite{BerGetVer92} and  \cite{Kos97}. Most of the definitions
  and properties of $\Cl{\g}$ and $Spin(\g)$ introduced here
  can be easily generalized to the Clifford algebra of any vector space
  with a non-degenerate symmetric bilinear form.

%%%%%%%%%%%%%%%%%%%%%%%%%%%%%%%%%%%%%%
%%% Action of Cl(g) on g %%%%%%%%%%%%%
%%%%%%%%%%%%%%%%%%%%%%%%%%%%%%%%%%%%%%

\subsection{The action of $\Cl{\g}$ on $\ext{\g}$}
  Using the ad-invariant inner product $B(\, , \, )$ on the Lie
  algebra $\g$, we define an $\Cl{\g}$-module structure on
  $\ext{\g}$ as follows. For $\forall x \in \g$, define
  \[
    \varrho(x)\cdot \alpha = x \wedge \alpha + \half \iota_{x} \alpha
               \quad \forall \alpha \in \ext{\g}
  \]

  where $\iota_x (y_1 \wedge \ldots \wedge y_k)
       = \sum_{i=1}^{k} (-1)^{i-1}B(x,y_i)\ y_1 \wedge \ldots
       \widehat{y_i} \ldots \wedge y_k $.
  It is easy to see that
  $\varrho(x_1)\varrho(x_2)+\varrho(x_2)\varrho(x_1)=B(x_1,x_2)$.
  So by the universal property of $\Cl{\g}$, $\varrho$ could be
  extended to the whole $\Cl{\g}$ by:
  \[
    \varrho(x_1 x_2 \ldots x_s)(\alpha)= \varrho(x_1) \circ \varrho(x_2)
        \ldots \circ \varrho(x_s)(\alpha).
  \]
  \begin{lemma}
    The inverse of quantization map $q: \ext{\g} \longrightarrow
    \Cl{\g}$ can be expressed in terms of $\varrho$ as $q^{-1}(x)=
    \varrho(x)\cdot 1$ for $\forall x \in \Cl{\g}$.
  \end{lemma}
  \begin{proof}
   Under an orthonormal basis ${e_1,\ldots, e_n}$ of $\g$, it is
   easy to see that, for any $1 \leq i_1 < \ldots < i_k \leq n$,
   \begin{align}
      \varrho(e_{i_1} \ldots e_{i_k})\cdot 1
           &=\varrho(e_{i_1}\ldots e_{i_{k-1}})\cdot e_{i_k} \notag \\
       \   &=\varrho(e_{i_1}\ldots e_{i_{k-2}})\cdot
              ( e_{i_{k-1}} \wedge e_{i_k} )  \notag  \\
       \   &=\cdots \cdots = e_{i_1} \wedge \ldots \wedge e_{i_k}. \notag
  \end{align}
   and conversely
   \[
      q(e_{i_1} \wedge \ldots \wedge e_{i_k}) = e_{i_1} \ldots
      e_{i_k}
   \]

   So the lemma follows from the fact that any element in $\Cl{\g}$ is a linear
   combination of $e_{i_1}\ldots e_{i_k}$ as above.
  \end{proof}

   People often call $q^{-1}$ \textit{symbol map}. It gives the
   isomorphism from the associated graded algebra $Gr(\Cl{\g})$ to
   $\ext{\g}$.

%%%%%%%%%%%%%%%%%%%%%%%%%%%%%%%%%%%%%%%%%%%%%%%%%%%%%%%%%%%%%%%
%%%%%% definition of Spin group and Spin Representation %%%%%%%
%%%%%%%%%%%%%%%%%%%%%%%%%%%%%%%%%%%%%%%%%%%%%%%%%%%%%%%%%%%%%%%

\subsection{Definition of Spin group and spin representation}
  We can think of $\Cl{\g}$ itself as a finite dimensional Lie
  algebra with the odd bracket defined by~\eqref{odd_bracket}. Notice the
  subspace $Cl^{(2)}(\g)$ is actually closed under the odd bracket, so
  $Cl^{(2)}(\g)$ is a Lie subalgebra of $\Cl{\g}$.
  In addition, let $SO(\g)$ be the special orthogonal group with
  respect to the inner product of $\g$. Its Lie algebra is denoted
  by $\mathfrak{so}(\g)$.
 \begin{lemma}
  $Cl^{(2)}(\g)$ is isomorphic to the Lie algebra
  $\mathfrak{so}(\g)$.
 \end{lemma}
 \begin{proof}
  We can identify $\g$ as $Cl^{(1)}(\g)$ and
  let $\tau_{\g} : Cl^{(2)}(\g) \longrightarrow \mathfrak{so}(\g)$ be
  defined by:
  \[
    \tau_{\g}(a)\cdot v = [a \, ,x]_{odd}, \;\, \forall\ a \in Cl^{(2)}(\g),
                                           x \in \g
  \]
   It is easy to check $\tau(a)$ does preserve $Cl^{(1)}(\g)$ and
   so defines a Lie algebra homomorphism from $Cl^{(2)}(\g)$ to
   $\mathfrak{gl}(\g)$. To see $\tau_{\g}(a)$ is in
   $\mathfrak{so}(\g)$, observe that
   \begin{align}
     B(\tau_{\g}(a)\cdot x,y) + B(x,\tau_{\g}(a)\cdot y)
       &= [\ [a,x]_{odd} ,y]_{odd} + [\ x,[a,y]_{odd}]_{odd} \notag \\
     \ &= -[\ a,[x,y]_{odd}]_{odd}= \ 0   \notag
   \end{align}
   The second equality is a consequence of the Jacobi identity in $(\Cl{\g}, [\; ,
   \,]_{odd})$.

     The map $\tau_{\g}$ must be an isomorphism, since it is
     injective by the non-degeneracy of the metric
     and since the dimensions of $Cl^{(2)}(\g)$ and
     $\mathfrak{so}(\g)$ are the same, namely $n(n-1)/2$.
 \end{proof}

    Using the map $\tau_{\g}$, any skew-symmetric matrix
    $A=(a_{ij}) \in \mathfrak{so}(\g)$ under an orthonormal
    basis $e_1, \cdots, e_n$ corresponds to the Cliiford element
   \begin{align}
     \tau^{-1}_{\g}(A) &= q( \half \sum_{i,j} (Ae_{i}) \wedge e_{j} )
                        = q( \half \sum_{i,j}a_{ji}e_{i} \wedge e_{j}) \\ \notag
                       &= q( -\half \sum_{i,j}a_{ij}e_{i} \wedge e_{j} )
                        = - \sum_{i<j}a_{ij}\xi_{i} \xi_{j} \notag
   \end{align}

\begin{definition}
  The group $Spin(\g)$ is the group obtained by exponentiating the
  Lie algebra $Cl^{(2)}(\g)$ inside the Clifford algebra $Cl(\g)$.
  The restriction of $\varrho$ on $Spin(\g)$ is called the \textit{Spin
  representation} of $Spin(\g)$.
\end{definition}

  The action $\tau_{\g}$ of $Cl^{(2)}(\g)$ on $\g$
  exponentiates to an orthogonal action still denoted by $\tau_{\g}$.
  So we have the following diagram.
\begin{equation}
  \begin{diagram}
   \node{\text{$Cl^{(2)}(\g)$}}
          \arrow[1]{e,t}{\text{exp}}
          \arrow{s,l}{\text{$\tau_{\g}$}}
   \node[1]{\text{$Spin(\g)\subset Cl(\g)$}}
          \arrow{s,r}{\text{$\tau_{\g}$}}
          \arrow{e,t}{\text{$\varrho$}}
   \node[1]{\text{End$(\wedge(\g))$}}  \\
   \node{\text{$\mathfrak{so}(\g)$}}
          \arrow[1]{e,t}{\text{exp}}
   \node[1]{\text{$SO(\g)$}}
  \end{diagram}
 \end{equation}

  It is well known that $\tau_{\g} : Spin(\g) \longrightarrow SO(\g)$ is a
  double covering map. But
  $Spin(\g)$ is much harder to handle than $SO(\g)$. For $SO(\g)$, we have a
  very nice representation in terms of matrices. We can
  investigate the structures of $SO(\g)$ using all kinds of
  decompositions of matrices. However, it is not easy to see how to factor
  the spin representation of an arbitrary element in $Spin(\g)$.
  In \cite{AM02}, A.Alekseev and E.Meinkenren
  constructed a very special factorization of $Spin(\g)$ which is interesting in
  many senses. We will see that the proof of the main theorem has to use
  this factorization.

  \begin{remark}
   The discussions above can be applied to any vector space $V$ with
   a non-degenerate symmetric bilinear form.
  \end{remark}

  \ \

%%%%%%%%%%%%%%%%%%%%%%%%%%%%%%%%%%%%%%%%%%%%%%%%%%
%%%%%% Factorization of spin representation %%%%%%
%%%%%%%%%%%%%%%%%%%%%%%%%%%%%%%%%%%%%%%%%%%%%%%%%%

  \subsection{Factorization of spin representation of $Spin(\g)$}
   For any $x \in Spin(\g)$, the spin representation $\varrho(x)$
   of x is a linear transformation on the vector space
   $\wedge(\g)$. In this section, we will write $\varrho(x)$ as
   the product of two special types of linear transformations on $\wedge(\g)$.
   We call this \textit{factorization of spin representation}.
   The two special transformations are:

   $ (1) \ \wedge(\g)  \longrightarrow \wedge(\g),  \;\;\; \alpha  \mapsto \beta \wedge
   \alpha \;\;\;\;\; $ for a fixed $\beta \in \wedge(\g)$

   $ (2) \ \wedge(\g)  \longrightarrow \wedge(\g), \;\;\; \alpha  \mapsto
         \iota_{\gamma}(\alpha) \;\;\;\;$ for a fixed $\gamma \in Cl(\g)$

   where $\iota_{x_1 \ldots x_s}(\alpha) = \iota_{x_1} \circ
   \ldots \circ \iota_{x_s}(\alpha)$.

  It is not clear how to do the factorization directly. So let us first
  introduce an auxiliary space.
  Since $\g$ is a finite-dimensional real vector space, the direct
  sum $W=\g \oplus \co{\g}$ carries a natural non-degenerate bilinear form
  \begin{align}
      B_{W}(x,y)          &= 0 , \;\, \forall\  x,y \in \g \ , \notag \\
      B_{W}(\alpha,\beta) &= 0 , \;\,  \forall\  \alpha , \beta \in
                               \co{\g} \ , \notag \\
      B_{W}(x,\alpha)     &= 2\alpha(x), \;\, \forall\  x\in \g, \alpha
                              \in \co{\g}
  \end{align}

  Let $Cl(W)$ be the Clifford algebra of $(W,B_W)$. The above discussion of $\Cl{\g}$
  can be applied to $Cl(W)$ without any change. So we have the Lie algebra
  isomorphism $\tau_{W}: Cl^{(2)}(W) \rightarrow \mathfrak{so}(W)$
  and the double cover of $SO(W)$ by $Spin(W)$.

   We can define an algebra representation $\pi$ of $Cl(W)$ on
   $\wedge(\g)$ by:
  \[
       \pi : Cl(W) \longrightarrow \mathfrak{gl}(\wedge(\g))
  \]
   Where generators $x\in \g$ act by wedge product and generators
   $\alpha \in \co{\g}$ act by contraction (denoted again by
   $\iota_{\alpha}$).

  Next, Let $\bar{\g}$ be the same vector space as $\g$ with metric
  $\bar{B}(\,\, ,\, )=-B(\,\, ,\, )$. Then we can construct an
  linear map $\kappa$ between $\g \oplus \bar{\g}$ and $W$ as follows.
  \[
     \kappa : \g \oplus \bar{\g} \rightarrow W,\: (x,y) \mapsto
             \left( x+y, B(\half(x-y),\:  -) \, \right)
  \]
  Assume the metric on $\g \oplus \bar{\g}$ is just the direct
  sum of their metrics. We have the following lemma.
 \begin{lemma}
  $\kappa$ is an isometry between $\g \oplus \bar{\g}$ and $W$ with
  respect to their metrics, its inverse is
  $\kappa^{-1}(x,B(y,-))=(\half x+y,\half x-y)$.
 \end{lemma}

  %\begin{remark}
   %$\g$ and $\bar{\g}$ can be identified with the real and imaginary
   %part of the complexfication $\g\otimes \C$,\ i.e.
   %$\g \otimes \C = \g \oplus i\g = \g \oplus \bar{\g}$.
   %In addition, if we identify $\co{\g}$ with $\g$ using the metric
   %$B(\: ,\,)$,
  %\end{remark}

  Using $\kappa$, we can identify $\g$ with a subspace of $W$.
  So $\Cl{\g}$ can be thought of as a subalgebra of $Cl(W)$.
  Correspondingly,
  $SO(\g)$ and $Spin(\g)$ are subgroups of $SO(W)$ and $Spin(W)$
  respectively.  Notice the restriction of $\kappa$ to $\g$ is:
   \[
     \kappa|_{\g}: \g \longrightarrow W,\:  x \mapsto \left(x,B(\half
     x,\;-)\right)
   \]
   So the representation $\varrho:Cl(\g) \rightarrow \mathfrak{gl}
   (\wedge(\g))$ is simply the restriction of
   $\pi : Cl(W) \longrightarrow \mathfrak{gl}(\wedge(\g))$ to $\Cl{g}$.

  The inclusion $h: SO(\g) \rightarrow SO(W)$ is given by:
 \[
    h: SO(\g) \rightarrow SO(W),\;
    C \mapsto \kappa\circ \left(
                               \begin{array}{cc}
                                  C & 0  \\
                                  0 & I
                                \end{array}
                         \right)
             \circ \kappa^{-1}
   =\left(
           \begin{array}{cc}
              \frac{1}{2} (C+I) & {C-I}  \\
              \frac{1}{4} (C-I) & \frac{1}{2}(C+I)
           \end{array}
     \right).
  \]

  Let $C=\exp(tA)$ where $A \in \mathfrak{so}(\g)$ and differential at $t=0$,
  we get the inclusion $dh: \mathfrak{so}(\g) \rightarrow
  \mathfrak{so}(W)$.

 \[
    dh: \mathfrak{so}(\g) \rightarrow \mathfrak{so}(W),\;
    A \mapsto \kappa\circ \left(
                               \begin{array}{cc}
                                  A & 0  \\
                                  0 & 0
                                \end{array}
                         \right)
             \circ \kappa^{-1}
   =\left(
           \begin{array}{cc}
              \frac{1}{2} A & A  \\
              \frac{1}{4} A & \frac{1}{2} A
           \end{array}
     \right).
  \]

  Then we can write down how $Cl^{(2)}(\g)$ and $Spin(\g)$ are included in
  $Cl^{(2)}(W)$ and $Spin(W)$. However, this won't directly give any
  decomposition of the spin representation of $Spin(\g)$.

  In \cite{AM02}, A.Alekseev and E. Meinrenken constructed a very
  interesting factorization of the spin representation of
  $Spin(\g)$. The idea is to
  first factor $SO(\g)$ inside $GL(W)$ and then lift the factorization
  to its double cover $Spin(\g)$ via the Lie algebra isomorphism
  $\tau_{W}: \mathfrak{so}(W) \rightarrow Cl^{(2)}(W)$. Please see
  the original paper for more details. Here, I just list the
  results they proved in \cite{AM02} for the future use. In fact,
  the statements in \cite{AM02} make sense for any vector space with
  an inner product.

  % Factorization of SO(g) on W
  \begin{thm}
  (Proposition 5.1 in \cite{AM02})
       Let $C \in SO(\g)$ with $\det(C-I) \neq 0$, and suppose that
    $D\in \mathfrak{so}(\g)$ is invertible and commutes with $C$. Then
    there is a unique factorization
   \begin{equation}\label{eq:factorization}
   h(C)=
   \left(\begin{array}{cc} I&0\\ E_1&I
         \end{array}  \right)
   \left(\begin{array}{cc} I&D\\0&I
         \end{array}  \right)
   \left(\begin{array}{cc} I&0\\ E_2&I
         \end{array}  \right)
   \left(\begin{array}{cc} R&0\\0&(R^{-1})^t
         \end{array}  \right)
   \end{equation}

 such that $E_1,E_2 \in \mathfrak{so}(\g)$ and $R \in GL(\g)$
 commute with $C$ and $D$. One finds
 \[
  E_1=\frac{1}{2} \frac{C+I}{C-I} - \frac{1}{D},\ \
  E_2=\frac{1}{D^2}\left(\frac{C-C^{-1}}{2}-D\right),\ \
  R=\frac{D}{I-C^{-1}}.
 \]
  \end{thm}

  \begin{remark}
    Notice $D \in \mathfrak{so}(\g)$ is invertible forces the
    dimension of $\g$ to be even. When the dimension of $\g$ is odd,
    we can define $W$ to be $\g \oplus \R
    \oplus \co{\g}$ instead and let the copy of $\R$ act on $\wedge(\g)$
    by scalar (see \cite{AM02} and \cite{GoodNolan98} for details).
    This won't effect our formula in~\eqref{eq:factorization} anyway.
  \end{remark}

 %%%%%%%%%%%%%%%%%%%%%%%%%%%%%%%%%%%%%%%%%
 %%%%%%%%% lifting SO(W) to Spin(W)  %%%%%
 %%%%%%%%%%%%%%%%%%%%%%%%%%%%%%%%%%%%%%%%%

  We know how to lift each component in the factorization ~\eqref{eq:factorization}
  to $Spin(W)$. In fact, under a basis $e_1,\cdots,e_n,e^1,\cdots,e^n$ of
  $W$, where $\{e_a\}$ is an orthonormal basis with respect to the
  metric $B(\, ,\,)$ on $\g$ and $\{e^a\}$ is its dual basis in
  $\co{\g}$, the lifting rules are:

  (1) For a matrix $D=(D_{ij})$ which represents a skew-adjoint linear map from
      $\co{\g}$ to $\g$, we have
     \[ \left(
          \begin{array}{cc}
             I & D \\
             0 & I
           \end{array}
      \right)
       \; \,
       \overset{\text{lift}}{\longrightarrow}
       \; \,
        \exp(-\half \sum_{i,j}D_{ij} e_i e_j) \; \in Spin(W).
     \]

  (2) For a matrix $E=(E_{ij})$ which represents a skew-adjoint linear map from $\g$ to
      $\co{\g}$, we have
     \[
      \left(
            \begin{array}{cc}
               I   &  0 \\
               E   &  I
         \end{array}
      \right)
      \; \,
       \overset{\text{lift}}{\longrightarrow}
       \; \,
       \exp(-\half \sum_{i,j}E_{ij}e^i e^j) \; \in Spin(W).
     \]

  (3) For any matrix $R \in GL(\g)$, there is a natural inclusion
     \[
        GL(\g) \rightarrow SO(W), \ \
      R\mapsto \left( \!
                      \begin{array}{cc}
                             R & 0  \\
                             0 & (R^{-1})^*
                      \end{array}\!\!\!
               \right).
     \]
     The lifting of this matrix to $Spin(W)$ is the same as lifting it to
     the \textit{metaplectic group} $Mp(2n,\R)$ (See \cite{Fol89} for the definition).
     It is not easy to write down the lifting explicitly in this case.
     But we do know how its
     lifting acts on $\wedge(\g)$ defined by $\pi$, which is enough to do our
     job of factoring spin representation of $Spin(\g)$ here.
     Let $\hat{R}$ be an element in $Mp(n,\g)$ (considered as a subgroup
     of $Spin(W)$) covering $R\in GL(\g)$. Then its action
     on $\wedge(\g)$ is given by:
      \begin{equation}\label{eq:R}
        \pi(\hat{R}) \cdot \alpha =\frac{R \cdot \alpha}{{|\det|}^{1/2}(\hat{R})}.
      \end{equation}
     where ${|\det|}^{1/2}: \, Mp(n,\R) \rightarrow \R^{\times}$ is a suitable
    choice of square root of $|\det|: \, GL(n,\g) \rightarrow \R^{+}$, and
    $R \cdot \alpha$ is defined by the unique extension of $R \in GL(\g)$ to an
    algebra automorphism of $\wedge(\g)$.

     \ \

   By the definition of $\pi$, for $\forall\ \alpha \in \wedge(\g)$,
   \begin{align}
      \pi(\exp(-\half \sum_{i,j}D_{ij} e_i e_j))\cdot \alpha &=
      \exp(-\half \sum_{i,j}D_{ij} e_i \wedge e_j) \wedge \alpha  \notag \\
      \pi(\exp(-\half \sum_{i,j}E_{ij} e^i e^j))\cdot \alpha &=
      exp(-\half \sum_{i,j}\iota_{E_{ij}e^i e^j})(\alpha)
   \end{align}

   So let $\lambda(D)=-\half \sum_{i,j}D_{ij}e_i \wedge e_j$ and
    $\gamma(E)=-\half \sum_{i,j}E_{ij}e^i e^j$, we get the
    following theorem.

 %%% Factorization of Spin(g)

 \begin{thm} \label{factorization}
 (Proposition 5.2 in \cite{AM02})
   Suppose $\hat{C}\in \Spin(\g)$ maps to $C \in SO(\g)$ with
   $\det(C-I)\neq 0$. Assume $h(C)$ has the factorization written
   in~\eqref{eq:factorization}. Then the operator $\varrho(\hat{C})$
   on $\wedge(\g)$ has the following factorization:
   \begin{equation}\label{eq:fac}
   \pi(\hat{C}) \cdot \alpha = \varrho(\hat{C}) \cdot \alpha =
      \frac{ \exp(\iota_{\gamma(E_1)}) \exp(\lambda(D)) \exp(\iota_{\gamma(E_2)})
             R \cdot \alpha }
           { {|\det|}^{1/2}(\hat{R}) }.
  \end{equation}
   In particular, when $\alpha =1$, we get
  \begin{equation}\label{eq:symb}
  \varrho(\hat{C}) \cdot 1= q^{-1}(\hat{C}) =
  \frac{\exp(\iota_{\gamma(E_1)})}
       { {|\det|}^{1/2}(\hat{R}) } \cdot \exp(\lambda(D)).
  \end{equation} \newline
 \end{thm}

 \begin{remark} \label{Explain}
 \end{remark}
 \begin{enumerate}
   \item This factorization actually makes sense for any vector space with
         an inner product.

   \item There are two ways to think of this decomposition. First,
         we can identify $\wedge(\g)$ with (odd) functions over $\co{\g}$
         and understand the theorem as giving a decomposition
         of some differential operators on the odd vector space $\co{\g}$.
         The second viewpoint is to think of $\wedge{\g}$ as the space of
         (odd) distributions supported at the origin of $\g$
         (see chapter 2). Correspondingly, elements in the Clifford
         algebra $Cl(\g)$ are thought of as the distributions with support at
         the identity element on a (super) Lie group of the
         (super) Lie algebra $(\g, [\;\, ,\;]_{odd})$.
         So ~\eqref{eq:symb} just tell us what the action of
         $\frac{\exp(\iota_{\gamma(E_1)})} { {|\det|}^{1/2}(\hat{R}) }$
         on the distribution $\exp(\lambda(D))$ is under the quantization map
         (which can be thought of as pushing forward of the distribution
         by the exponential map).

   \item We will see that if we set $C=\exp(adx)$ and $D=adx$, the theorem
         gives the expression of our quantization map $\mathcal{Q}$.
  \end{enumerate}

   %%%%%%%%%%%%%%%%%%%%%%%%%%%%%%%%%%%
   %%% Proof of the Main Theorem %%%%%
   %%%%%%%%%%%%%%%%%%%%%%%%%%%%%%%%%%%

 \subsection{Proof of the Main theorem}
    The basic idea is to apply the factorization of spin
   representation we established in the last
   chapter to some special elements in $Spin(\g)$ and get a highly
   nontrivial relation in the Clifford algebra. Now Let us begin.

  \begin{proof}
     \ \
    Recall in chapter 2,
    we used $v^a$ and
    $\theta^a$ to denote the generators of $S(\co{\g})$ and $\wedge{\co{\g}}$
    in $W_{\g}$. And we use $\widehat{v^a}$ and
    $\theta^a$ to denote the generators of $S(\co{\g})$ and $\wedge{\co{\g}}$
    in $W^{K}_{\g}$. Let their dual generators in $S(g)$ and $\wedge(\g)$
    be $e_a,\ \widehat{e_a}$ and $\bar{e_{a}}$. By definition of
    $v^a$ and $\widehat{v^a}$, we have $\widehat{e_a}=e_a-\half
    f_{abc}\bar{e_b}\bar{e_c}$. So
    $e_a$ commutes with $\bar{e_{b}}$ while $\widehat{e_a}$ doesn't.
    Let $\g_{0}$ and $\g_{1}$ be the Lie algebras
    $\widetilde{T\g[1]}^{even}$ and $\widetilde{T\g[1]}^{odd}$
    respectively.

   In formula~\eqref{eq:fac}, for $\mu= \mu^a \widehat{e_a} \in \g_{0}$,
   set $C=\exp_{\mathfrak{so}(\g)}(ad\mu), \hat{C}=\exp_{Cl(\g)}(\half
   \sum_{i,j}(ad\mu)_{ij}\bar{e_i}\bar{e_j})$
   , and $D=ad\mu$. Then we get:
   \[
     \frac{1}{{|\det|}^{1/2}(\hat{R})} = j^{\half}(ad\mu), \;\:
     E_1= f(ad\mu),\;\
   \]

   \[
           f(s)= \frac{1}{2}\frac{e^s +1}{e^s -1} -
                      \frac{1}{s} = \frac{dlnj(s)}{ds}
   \]
    \

   Then $\exp(\iota_{\gamma(E_1)})$ is exactly the
   $exp \left( \half T_{ab}\iota_{a}\iota_{b} \right)$ in the
   quantization map $\mathcal{Q}$. If we
   work in the Fourier transformed picture of the usual interpretation of Duflo
   map, elements in $\g_0$ and $\g_1$ are considered as even and odd
   distributions supported at the origin respectively.

   Let $\iota_{\mathcal{C}(x)}$ denote multiplying the function
   $j^{\half}(x) \cdot exp \left( \half T_{ab}(x) \bar{e_a}\bar{e_b} \right)$ to a
   distribution, then the factorization~\eqref{eq:symb} reads:
   \begin{equation} \label{eq:specialfac}
      q^{-1} \left( \exp_{Cl(\g)}(-\half \sum_{i,j}(ad\mu)_{ij} \bar{e_i}\bar{e_j})
                   \right)
      = \iota_{\mathcal{C}(\mu)} \cdot \exp_{\wedge(\g)}\left(-\half (ad\mu)_{ij} \bar{e_i}
      \bar{e_j}\right)    \end{equation}

    (The $\mu$ part serve as even functions). Furthermore, the
    main theorem in \cite{AM02} says the following is also true.
    \begin{align}
       & \quad\:\: q^{-1} \left(
                   \exp_{Cl(\g)} \left(
                                        -\half \sum_{i,j}((ad\mu)_{ij} \bar{e_i}\bar{e_j})
                                        + \nu^{a}\bar{e_a}
                                 \right)
                        \right)   \notag \\
       &= \iota_{\mathcal{C}(\mu)} \circ \exp_{\wedge(\g)}
                             \left(
                                 -\half \sum_{i,j}((ad\mu)_{ij} \bar{e_i}\bar{e_j})
                                  + \nu^{a}\bar{e_a}
                             \right)  \notag
    \end{align}
     The proof of this is essentially an application of
     theorem~\eqref{factorization} to a slightly larger space.(See \cite{AM02} for detail).
     In addition, since
     \[
         \sum_{i,j}(ad\mu)_{ij} \bar{e_i}\bar{e_j}= \sum_{a,b,c}  \mu^{a}
        f_{abc}\bar{e_b}\bar{e_c}.
     \]
     We get
     \begin{align} \label{eq:strongfac}
       & \quad\:\:  q^{-1} \left(
                   \exp_{Cl(\g)} \left(
                                        -\half \sum_{a,b,c}(\mu^{a} f_{abc}\bar{e_b}\bar{e_c})
                                        + \nu^{a}\bar{e_a}
                                 \right)
                           \right)  \notag \\
       &= \iota_{\mathcal{C}(\mu)} \circ \exp_{\wedge(\g)}
                             \left(
                                 -\half \sum_{a,b,c}(\mu^{a} f_{abc}\bar{e_b}\bar{e_c})
                                  + \nu^{a}\bar{e_a}
                             \right)
    \end{align}

     Since we will deal with several different algebras in our proof here, It is convenient
    to just use some subscripts under $exp$ to indicate the multiplications in
    different algebras and omit the $\wedge$ symbols between odd elements.

    First, Let us show that the theorem is true for the special element
   $ \exp_{S(\widetilde{T\g[1]})}(\sum_{a}(\mu^a \widehat{e_a} + \nu^a \bar{e_a}))
   $ where $\mu_a$ and $\nu_a$ are parameters, i.e. we have to show :
   \begin{align} \label{case1}
     & \quad\: \mathcal{Q} \circ \sigma_{0} \left( \exp_{S(\widetilde{T\g[1]})}
     (\sum_{a}(\mu^a \widehat{e_a} +
   \nu^a \bar{e_a})) \right)  \notag \\
     &= \sigma_1 \circ Duflo \left (\exp_{S(\widetilde{T\g[1]})}
    (\sum_{a}(\mu^a \widehat{e_a} + \nu^a \bar{e_a})) \right)
   \end{align}

    Since Duflo map here is just the super-symmetrization, so the right side
    of~\eqref{case1} is:
   \begin{align}
       & \;\;\;\;\; \sigma_1  \circ Duflo \left( \exp_{S(\widetilde{T\g[1]})} \sum_{a}
          (\mu^a \widehat{e_a} + \nu^{a} \bar{e_a} ) \right)  \notag \\
       &= \sigma_1 \left(\exp_{U(\widetilde{T\g[1]})} \left( \sum_{a}
          (\mu^a \widehat{e_a} + \nu^a \bar{e_a}) \right) \right)  \notag \\
       &= \exp_{U(\widetilde{T\g[1]})} \left( \sum_{a} \left( \mu^a (e_a - \half f_{abc} \bar{e_b}\bar{e_c})
           + \nu^a \bar{e_a} \right) \right)  \notag \\
       &= \exp_{U(\widetilde{T\g[1]})} \left( \sum_{a} \left(
              \mu^{a} e_a + ( -\half \mu^{a} f_{abc}\bar{e_b}\bar{e_c}
           + \nu^a \bar{e_a}) \right) \right) \notag \\
      &=  \exp_{U_{\g}} \left( \sum_{a}  \mu^{a} e_a  \right) \otimes\,
          \exp_{Cl(\g)} \left( \sum_{a} (-\half \mu^{a} f_{abc}\bar{e_b}\bar{e_c}
           + \nu^a \bar{e_a} ) \right).  \notag
   \end{align}

    The last equality is because $e_a$ commutes with
    $\bar{e_a}$.

    Notice the quantization map $\mathcal{Q}=(\chi \otimes q) \circ
    \iota_{\mathcal{C}(x)}$, the left side of~\eqref{case1} is:
    \begin{align}
     & \;\;\;\;\;\;\; \mathcal{Q} \circ \sigma_{0} \left(
         \exp_{S(\widetilde{T\g[1]})} \sum_{a}(\mu^a \widehat{e_a} + \nu^{a} \bar{e_a} )
                                  \right)  \notag \\
     & \;\; = \mathcal{Q} \cdot \exp_{S(\widetilde{T\g[1]})}
                       \left( \sum_{a}
                                 \left( \mu^a (e_a - \half f_{abc} \bar{e_b}\bar{e_c})
                                    + \nu^a \bar{e_a}  \right)
                       \right)  \notag \\
     & \;\;= (\chi \otimes q) \circ \iota_{\mathcal{C}(x)}
                              \left(  \exp_{S(\g)} \left( \sum_{a} \mu^a e_a \right)
                                      \otimes \, \exp_{\wedge(\g)}
                                 \left(
                                      \sum_{a} ( -\half \mu^{a} f_{abc}\bar{e_b}\bar{e_c} + \nu^{a}
                                                \bar{e_a} )
                                 \right)
                              \right) \notag \\
     & \;\; = \chi \left(
                   \exp_{S(\g)}
                      \left(
                               \sum_{a} \mu^a e_a
                      \right)
              \right) \otimes\:
        q \left( \iota_{\mathcal{C}(\mu)} \cdot  \exp_{\wedge(\g)}
                \left(
                     \sum_{a} (- \half \mu^{a} f_{abc}\bar{e_b}\bar{e_c} + \nu^{a}
                       \bar{e_a} )
                 \right)
          \right)                         \notag \\
    & \overset{~\eqref{eq:strongfac}}{=} \exp_{U_{\g}}
                                 \left(
                                      \sum_{a}  \mu^{a} e_a
                                 \right)  \otimes \,
          \exp_{Cl(\g)} \left( \sum_{a} (-\half \mu^{a} f_{abc}\bar{e_b}\bar{e_c}
           + \nu^a \bar{e_a} ) \right). \notag
    \end{align}

    Here, we think of elements in $\sym{\widetilde{T\g[1]}}$
    as distributions supported at origin (see chapter 2).
    Hence we have proved the equation~\eqref{case1}.

    By comparing the coefficients of the parameters $\mu_a$ and
    $\nu_a$ of both sides in~\eqref{case1},
    we prove the theorem for any elements in $\sym{\widetilde{T\g[1]}}$.
    \newline
 \end{proof}

 %%%% end of Chapter 4 %%%%%%%%
 %Proof of the Main Theorem

\large
\section{\textbf{Jacobi diagrams and diagrammatic proof of Duflo isomorphism }}
\normalsize

This chapter is an explanation of the work \cite{BTT03} of
Bar-Natan, Le and Dylan Thurston.

%%%%%%%%%%%%%%%%%%%%%%%%%%%%%%%%%%%%%%%%%%%%%%
%%%%%%%% Tensors and Jacobi diagrams %%%%%%%%%
%%%%%%%%%%%%%%%%%%%%%%%%%%%%%%%%%%%%%%%%%%%%%%

\subsection{Tensors and Jacobi diagrams}
 In general, a tensor with $n$ free indices will be represented by a graph with $n$
 legs. The indices of a tensor can belong to different vector spaces or their
 dual spaces. Correspondingly, the legs of the graph should be labeled to indicate
 the vector spaces and distinguish a vector space and its dual. For example, a
 matrix $A \in Hom(V,V)= V^* \otimes V$ can be represented by
 the graph in figure~\eqref{fig_matrix}. By convention, the data flow in the direction
 of the arrows, so the incoming arrow is the $V^*$ factor and the outgoing arrow is
 the $V$ factor. In this paper, we mainly deal with the tensors over a Lie algebra $\g$
 with an invariant metric  $B(\, ,\, )$.

 \begin{figure}
  \begin{equation*}
   \vcenter{
            \hbox{
                  \mbox{$\input{matrix.pstex_t}$}
                 }
           } \quad \quad \longleftrightarrow \quad A \in V^* \otimes V
  \end{equation*}
  \caption{Pictorial representation of A  \label{fig_matrix}}
 \end{figure}

\begin{figure}
  \begin{equation*}
   \vcenter{
            \hbox{
                  \mbox{$\input{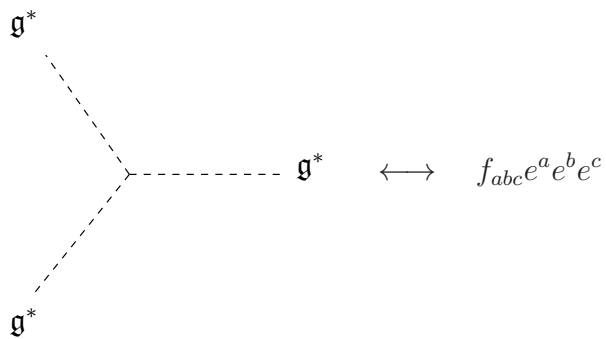}$}
                 }
           } \quad \quad \longleftrightarrow \quad f_{abc}e^{a}e^{b}e^{c}
  \end{equation*}
  \caption{Pictorial representation of Lie bracket \label{fig_3_fork}}
 \end{figure}

\begin{figure}
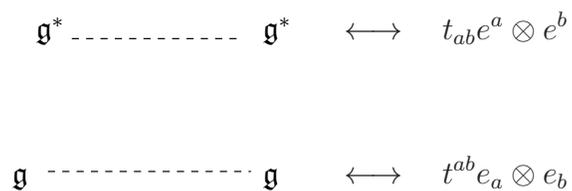

  \begin{align}
   \vcenter{
            \hbox{
                  \mbox{$\input{line_seg.pstex_t}$}
                 }
           } \quad \quad \longleftrightarrow \quad t_{ab}e^{a} \otimes e^{b} \notag \\
        \quad \quad \quad   \notag        \\
        \quad \quad \quad   \notag        \\
   \vcenter{
            \hbox{
                  \mbox{$\input{line_seg1.pstex_t}$}
                 }
           } \quad \quad \longleftrightarrow \quad t^{ab}e_{a} \otimes e_{b} \notag
  \end{align}
  \caption{Pictorial representation of the invariant metric B(\, ,\, ) and its dual \label{fig_line_seg}}
 \end{figure}

  Suppose $e_1, \ldots, e_n$ is a basis of $\g$ and $e^1,
  \ldots, e^n$ is the dual basis in $\co{\g}$.
  Suppose $[e_a,e_b]=f_{ab}^{c}e_c$ and $B(e_a,e_b)=t_{ab}$.
  The matrix $(t_{ab})$ is invertible.
  We use $(t^{ab})$ to denote its inverse matrix, which defines
  the dual metric of B(\, ,\, ) on $\co{\g}$. Then
  the invariant metric B(\, ,\, ) and its dual metric could be
  represented by diagrams in figure~\eqref{fig_line_seg}.
  It is easy to show that $f_{abc}=f_{ab}^{k}t_{kc}$ is totally
  anti-symmetric in its indices. We use a fork (see figure~\eqref{fig_3_fork})
  to represent it.

\begin{figure}
  \begin{align}
   \vcenter{
            \hbox{
                  \mbox{$\includegraphics{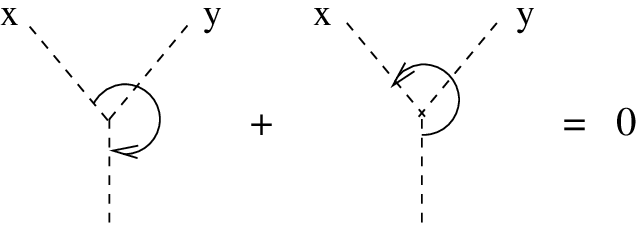}$}
                 }
           }   \notag \\
   \vcenter{
            \hbox{
                  \mbox{$\includegraphics{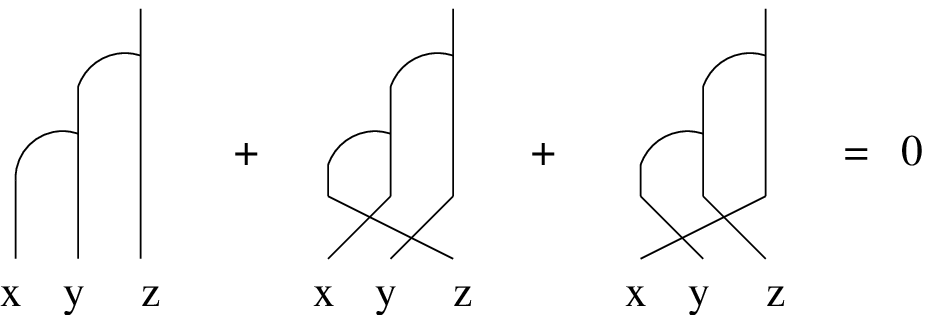}$}
                 }
           }   \notag
  \end{align}
  \caption{Pictorial representation of antisymmetry and Jacobi identity of Lie bracket
           \label{fig_Jacobi}}
 \end{figure}

The antisymmetry and Jacobi identity for the bracket of the Lie
algebra $\g$
\[
    [x,y]+[x,y]=0, \quad \quad
    [[x,y],z]+[[z,x],y]+[[y,z],x]=0 \quad \; \forall x,y,z \in \g
\]
  can be expressed graphically as in figure~\eqref{fig_Jacobi}.

  The two relations of diagrams in figure~\eqref{fig_Jacobi} are
  called \textit{antisymmetry} and IHX relations.

  Notice that the trivalent vertices in a diagram have to be
  oriented, otherwise it is not clear how to read the corresponding
  tensor from the diagram.

 \begin{definition}
   An open Jacobi diagram (also called uni-trivalent graph) is a
   vertex-oriented uni-trivalent graph, i.e., a graph with only
   univalent and trivalent vertices where each trivalent vertex is
   oriented. The univalent vertices are called legs. In
   planar pictures, the orientation on the edges incident on a
   vertex is always anticlockwise.
 \end{definition}

 \begin{definition} \label{Jacobi_Diagrams}
   Let $\mathcal{B}^f$ be the vector space spanned by Jacobi
   diagrams modulo the IHX relation and the antisymmetry relation.
   The degree of a diagram in $\mathcal{B}^f$ is half of the number
   of vertices (trivalent and univalent) of the diagram. Let
   $\mathcal{B}$ be the completion of $\mathcal{B}^f$ with respect
   to the degree.
 \end{definition}

  There are some remarkable relations between the space of
  Jacobi diagrams and the space of \textit{Vassiliev invariants}. Roughly
  speaking, any finite type weight system on $\mathcal{B}$ (i.e. a real value function
  on $\mathcal{B}$ that vanishes on diagrams with degree greater than
  a certain integer) corresponds to a
  Vassiliev invariant. The correspondence is established by M.Kontsevich using
  his famous Kontsevich integral. See Kontsevich's original paper \cite{Ko93}
  and Dror Bar-Natan's paper \cite{Bar95} for the complete exposition on this
  topic.

  If we identify $\co{\g}$ and $\g$ using the metric, Any Jacobi diagram
  $D$ gives a tensor over $\g$ in the following way. First we can
  decompose $D$ into several copies of forks and
  bars, ignoring the crossings between them. Each of the forks and bars
  is canonically associated to an invariant tensor as in
  figure~\eqref{fig_3_fork} and figure~\eqref{fig_line_seg}.
  If contracting the legs of these forks and bars, we get a tensor.
  It is easy to see the tensor we get does not depend on the way we
  decompose the diagram, hence we denote it by $T_D$.
  Next, since the symmetric algebra $\sym{\g}$ is a quotient of $T(\g)$,
  let $p: T(\g) \longrightarrow \sym{\g}$ be the quotient map. It is easy to see that
  $p(T_D) \in \sym{\g}$ will correspond to the diagram $D$ if we don't order the legs
  of $D$. That is to say every Jacobi diagram with unordered legs gives a symmetric
  tensor over $\g$. Unfortunately, not every element in $\sym{\g}$ can be represented
  in this way because of the following lemma.

\begin{lemma}
 The symmetric tensor associated to a diagram in $\mathcal{B}$ is
 always invariant under the adjoint action of $\g$.
\end{lemma}
\begin{proof}
 The tensor associated to a diagram in $\mathcal{B}$ is made by
 contracting copies of structure constants and invariant metric
 tensors of the Lie algebra, which are all invariant under $ad(\g)$.
\end{proof}

\begin{remark}
  Actually, not every invariant symmetric tensor can be represented
  by diagrams in $\mathcal{B}$, see \cite{Vogel} for examples.
\end{remark}

\begin{figure}
  \begin{equation*}
   \vcenter{
            \hbox{
                  \mbox{$\input{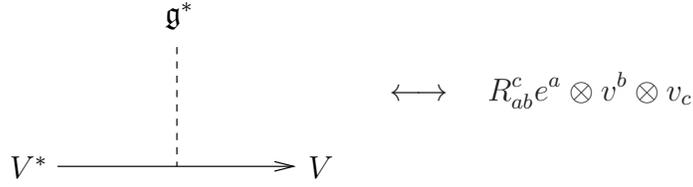}$}
                 }
           } \quad \quad \longleftrightarrow \quad R_{ab}^{c}e^{a} \otimes v^{b} \otimes v_{c}
  \end{equation*}
  \caption{Jacobi diagram of a representation $V$ of $\mathfrak{g}$\label{fig_representation}}
 \end{figure}

\begin{figure}
  \begin{equation*}
   \vcenter{
            \hbox{
                  \mbox{$\includegraphics{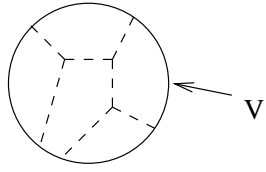}$}
                 }
           }
  \end{equation*}
  \caption{An example of Jacobi diagram \label{fig_Sg_example}}
 \end{figure}

 We can another ingredient into the Jacobi diagrams.
 If we have a representation $(V,\pi)$ of $\g$, it defines an
 element $R \in V^* \otimes \co{\g} \otimes V$, which is
 represented by figure~\eqref{fig_representation} where
 we use a different type of line to denote $V$. The fact that $V$
 is a representation of $\g$ imposes an IHX type relation in this
 setting. Now, we can draw diagrams with different type of edges,
 e.g. figure~\eqref{fig_Sg_example}, which allows us to construct more
 invariant tensors over $\g$.

 In addition, $\mathcal{B}$ becomes a commutative algebra if we
 define the product of two diagram to be the disjoint union $\sqcup$ of them.
 This corresponds exactly to the algebra structure on $\sym{\g}$.

%%%%%%%%%%%%%%%%%%%%%%%%%%%%%%%%%%%%%%%%%%%%%
%%%%%%%%% Based Jacobi diagram %%%%%%%%%%%%%%
%%%%%%%%%%%%%%%%%%%%%%%%%%%%%%%%%%%%%%%%%%%%%

\subsection{Based Jacobi diagrams}
\begin{figure}
  \begin{equation*}
   \vcenter{
            \hbox{
                  \mbox{$\includegraphics{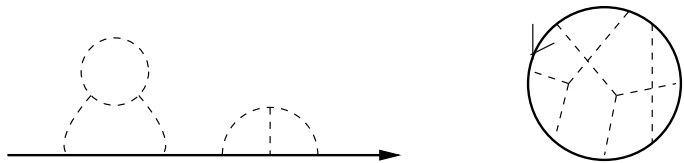}$}
                 }
           }
  \end{equation*}
  \caption{Examples of based Jacobi diagrams \label{fig_Ug_example}}
 \end{figure}

\begin{figure}
  \begin{equation*}
   \vcenter{
            \hbox{
                  \mbox{$\includegraphics{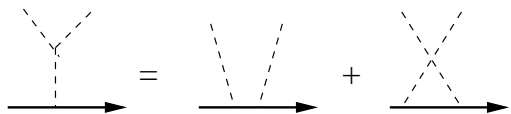}$}
                 }
           }
  \end{equation*}
  \caption{STU relation of based Jacobi diagrams \label{fig_STU}}
 \end{figure}

 If we order the legs of a Jacobi diagram $D$, it will represent
 a tensor with noncommutative indices. To remember the ordering, we
 can glue the legs of $D$ to a connected 1-manifold $X$
(a circle or a oriented line), see figure~\eqref{fig_Ug_example}.
 In this case, we can't commute any two legs that are attached to $X$
 without changing the diagram. If we impose another diagrammatic relation,
 called STU relation (see figure~\eqref{fig_STU}), between any two legs
 attached to $X$, then the diagram represents an element in
 the universal enveloping algebra $U(\g)=T(\g)/<xy-yx=[x,y]>$.
 We call this new type of diagrams  \textit{based Jacobi diagrams} on $X$.

 \begin{definition}
  For an oriented connected 1-manifold $X$, let $\mathcal{A}^f(X)$ be the
  vector space spanned by Jacobi diagrams based on $X$ modulo the
  antisymmetry, IHX and STU relations.  The \textit{degree} of a
  diagram in $\mathcal{A}^f(X)$ is half the number of its vertices.
  Define $\mathcal{A}(X)$ to be the completion of $\mathcal{A}^f(X)$
  with respect to the degree.
 \end{definition}

  If we have a representation $V$ of the Lie algebra $\g$, we can
  draw diagrams with dashed lines and solid lines with legs
  attached to $X$.

 \textbf{Warning}: the solid line denoting $V$ and the based 1-manifold
  are not the same thing. We don't put any representation of $\g$ on
  the based 1-manifold.

 The algebra of Jacobi diagrams based on a oriented circle
 $\mathcal{A}(\circlearrowleft)$ and Jacobi diagrams based
 on an oriented line $\mathcal{A}(\longrightarrow)$ are actually
 isomorphic as algebras (see \cite{Bar95}). So we can just use
 $\mathcal{A}$ to denote $\mathcal{A}(\circlearrowleft)$ or
 $\mathcal{A}(\longrightarrow)$.

\begin{figure}
  \begin{equation*}
   \vcenter{
            \hbox{
                  \mbox{$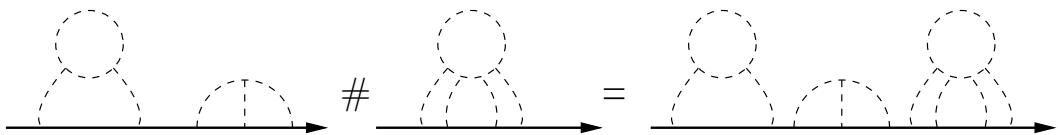$}
                 }
           }
  \end{equation*}
  \caption{Product of based Jacobi diagrams \label{fig_Ug_product}}
 \end{figure}

 Similar to $\mathcal{B}$, we can define an algebra structure
 on $\mathcal{A}:$ take two based Jacobi diagrams $D_1,D_2$
 and place the legs of
 $D_1$ before the legs of $D_2$ in the total ordering of legs. We denote
 it by $D_1$\#$D_2$.
 In diagrams, it is just connecting the base lines of $D_1$ and
 $D_2$ (see figure~\eqref{fig_Ug_product}).

\begin{figure}
  \begin{equation*}
   \vcenter{
            \hbox{
                  \mbox{$\includegraphics{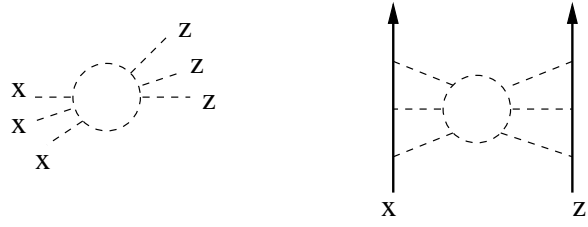}$}
                 }
           }
  \end{equation*}
  \caption{Colored Jacobi diagrams \label{fig_Colored_Diagrams}}
 \end{figure}

 Another useful type of Jacobi diagrams is the \textit{Jacobi diagram with colored
 legs} in which we give different colors to legs of a Jacobi diagram and don't
 distinguish legs with the same color. When we glue the colored legs to an
 oriented line or a circle, the line or the circle is also colored.
 See figure~\eqref{fig_Colored_Diagrams}. We use
 $\mathcal{A}(\ast_{x_1}, \ldots, \ast_{x_n})$ to denote the Jacobi diagrams with
 legs colored by $x_1, \ldots, x_n$.

 For more discussion on various kinds of Jacobi diagrams, see \cite{Bar95}
 and \cite{BTT03}.

%%%%%%%%%%%%%%%%%%%%%%%%%%%%%%%%%%%%%%%%%%%%%%%%%%%%%%%%%%
%%%%%%%%% Diagrammatic proof of Duflo isomorphism %%%%%%%%
%%%%%%%%%%%%%%%%%%%%%%%%%%%%%%%%%%%%%%%%%%%%%%%%%%%%%%%%%%

\subsection{Diagrammatic proof of Duflo isomorphism}
 In \cite{BTT03}, Dror Bar-Natan, Thang T. Q. Le and Dylan P.
 Thurston gave an interesting diagrammatic analogue of Duflo isomorphism in
 the world of Jacobi diagrams. Their proof essentially uses some
 special properties of the Kontsevich integral of unknot. For
 the definition and properties of the Kontsevich integral of knots
 and tangles, please see \cite{Bar95}, \cite{BTT03} and \cite{Ohto02}.

\begin{figure}
  \begin{equation*}
   \vcenter{
            \hbox{
                  \mbox{$\includegraphics{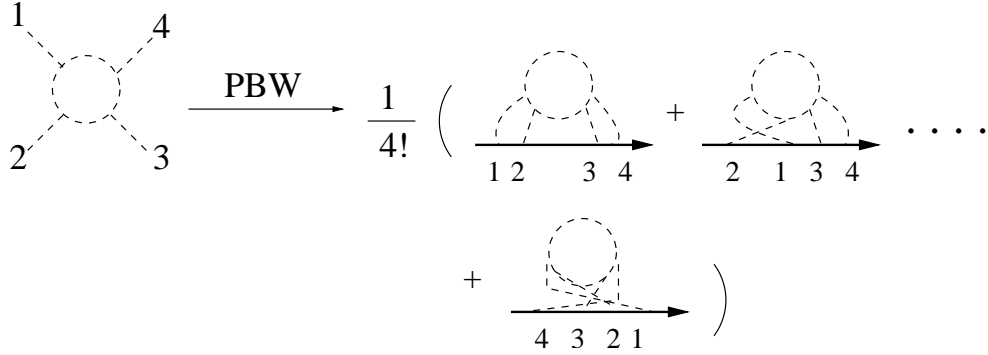}$}
                 }
           }
  \end{equation*}
  \caption{PBW isomorphism for Jacobi diagrams \label{fig_PBW}}
 \end{figure}

 We can interpret the Duflo isomorphism in terms of Jacobi diagrams as follows.
 First, the Poincar\'{e}-Birkhoff-Witt map can be described diagrammatically
 as averaging all ways of gluing all the legs of a Jacobi diagram to an oriented line
 (see figure~\eqref{fig_PBW}).

 \begin{definition}
  For diagrams $C, C' \in \mathcal{B}$ so that $C$ has no struts (components like
  \includegraphics[scale=0.8]{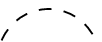}), the
  inner product of $C$ and $C'$ is defined by
  \begin{equation} \label{Diagram_Drivative}
    \langle C, C' \rangle =
    \begin{cases}
     \text{the sum of all ways of gluing all} &\text{\quad if $C$ and $C'$ have the same}\\
     \text{legs of C to all legs of D,}        &\text{\quad number of legs,} \\
     0                                        & \quad
     \text{otherwise}
    \end{cases}
  \end{equation}
   If $C$ and $C'$ are colored Jacobi diagrams, we require that only when the legs
   from $C$ and $C'$ have the same color, can they be glued together. So in this case,
   $C$ and $C'$ must have the same number of colored legs in each color to make their inner
   product non-zero.
 \end{definition}

  If we write $\iota_{C} \stackrel{def}{=} \langle C \, , \: \rangle
  $, the PBW map $\chi$ can be pictorially thought of as
  $\iota_{\Gamma}$, where $\Gamma$ is a colored Jacobi diagram
  in $\mathcal{A}(\ast_z,\ast_x)$ (see figure~\eqref{fig_Chi}).

 \begin{figure}
  \begin{equation*}
   \vcenter{
            \hbox{
                  \mbox{$\input{Chi_1.pstex_t}$}
                 }
           }
  \end{equation*}
  \caption{The diagram $\Gamma$ for PBW map $\chi$
  \label{fig_Chi}}
 \end{figure}

\begin{figure}
  \begin{equation*}
   \vcenter{
            \hbox{
                  \mbox{$\input{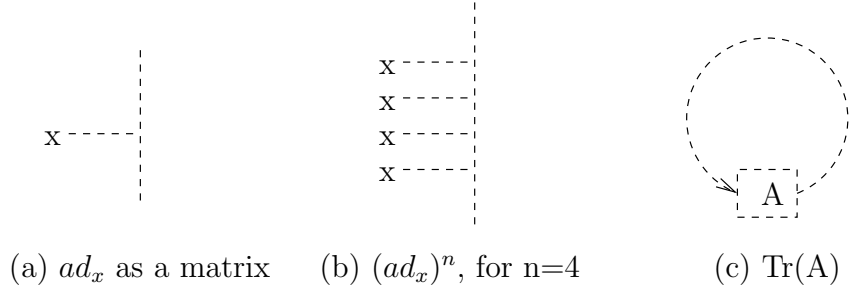}$}
                 }
           }
  \end{equation*}
  \caption{Building blocks of wheels \label{fig_ad_x}}
 \end{figure}

 Next, we interpret the $j^{\half}$ in Duflo isomorphism in the language
 of Jacobi diagrams. For any $x \in \g$,
 Think of $ad_x$ as a matrix, which is represented by the diagram
 in figure~\eqref{fig_ad_x}(a). Then $(ad_x)^n$ can be denoted by
 the diagram in figure~\eqref{fig_ad_x}(b). In addition, taking trace of a
 matrix is just connecting the input and the output, see
 figure~\eqref{fig_ad_x}(c). So we have the diagrammatic
 representation of $j^{\half}$ as following.
 \begin{align}
  j^{\half}(x) &= det^{\half} \left( \frac{\sinh(\half
                 ad_{x})}{\half ad_{x}} \right)
                = \exp \left(
                       \half tr
                       \left(
                             \log \frac{\sinh(\half ad_{x})}{\half ad_{x}}
                       \right)
                      \right) \notag  \\
               &= \exp \left(
                        \half tr \left(
                                  \sum_{n=0}^{\infty}
                                  b_{2n}(ad_x)^{2n}
                                 \right)
                       \right)  \notag \\
               &= \exp \left(
                         \frac{1}{48} \includegraphics[scale=0.5]{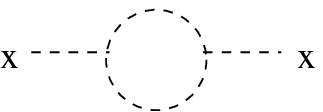}
                       - \frac{1}{5760} \includegraphics[scale=0.5]{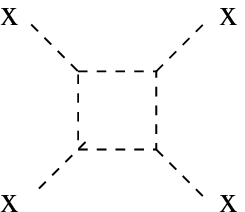}
                       + \frac{1}{362880} \includegraphics[scale=0.5]{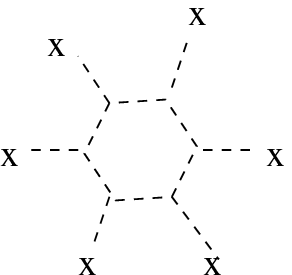}
                       - \ldots
                       \right)  \label{omega} \overset{\vartriangle}{=} \Omega_x
 \end{align}

  The diagram $\Omega_x$ is a very important Jacobi diagram. It is proven in
  \cite{BTT03} to be the Kontsevich integral of the unkot -- the only Kontsevich
  integral of a knot that people can calculate completely so far!
  We use $\Omega_x$ to denote this diagram (
  the subscript x here means the legs of the diagram is colored by x, so it is actually
  an element in $\mathcal{A}(\ast_x)$).

  \begin{remark}
   Notice the diagrammatic
   interpretation of the Duflo isomorphism in figure~\eqref{fig_Chi}
   and~\eqref{fig_ad_x} makes sense even without the metric on $\g$.
  \end{remark}

 In addition, Since $j^{\half}(x)$ acts as an differential operator on
 $\sym{\g}$ in
 the Duflo isomorphism, we have to define how differential operators are
 represented in the world of Jacobi diagrams.

 \begin{definition}
  For a $C \in \mathcal{B}$ without struts, the operation of
  applying $C$ as a differential operator,
  denoted by $\partial_{C}: \mathcal{B} \longrightarrow \mathcal{B}$,
  is defined to be
  \begin{equation} \label{DigramOperator}
   \partial_{C}=
   \begin{cases}
     0  &\text{if C has more legs than D,} \\
    \text{the sum of all ways of gluing all} &\text{otherwise.} \\
    \text{legs of C to some(or all) legs of D}
   \end{cases}
  \end{equation}

   If $C$ and $D$ are colored Jacobi diagrams, we can only glue legs of the
   same color in above definition. In addition, let $\varnothing$ denote
   the empty diagram, then $\partial_{\varnothing}(D) = D$.

 \end{definition}

  With this definition, we can interpret the action of $j^{\half}$ on symmetric
  algebra $\sym{\g}$ diagrammatically as $\partial_{\Omega}$ on $\mathcal{B}$.

  Before stating the diagrammatic analogue of Duflo isomorphism, we need to
  introduce two more operations in Jacobi diagrams.

 \begin{definition} \label{coloring}
 For any Jacobi diagram $C$ (no matter if it is colored), $(C)_x$ is just $C$
 with all its legs colored by x, ignoring the original coloring of $C$.
 $\Delta_{xy}C$ is defined to be the sum of all possible
 different colorings on the legs of $C$ by two colors $x$ and $y$, ignoring the
 original coloring of $C$.
 \end{definition}

 \begin{lemma} \label{duality}
  For any diagrams $C, D_1, D_2 \in \mathcal{B}$,
  \[
    \langle C, D_1 \sqcup D_2 \rangle = \langle \Delta_{xy}C, (D_1)_x
    \sqcup (D_2)_y  \rangle
  \]
 \end{lemma}
 \begin{proof}
  Obvious.
 \end{proof}

  we have the following theorem which is established in \cite{BTT03}.

 \begin{thm}(Wheeling)
   The map $\Phi = \partial_{\Gamma} \circ \partial_{\Omega_x} : \mathcal{B}
   \longrightarrow \mathcal{A}$ is an algebra homomorphism with respect to
   the algebraic structure in $\mathcal{B}$ and $\mathcal{A}$.
 \end{thm}

\begin{figure}
  \begin{equation*}
   \vcenter{
            \hbox{
                  \mbox{$\includegraphics{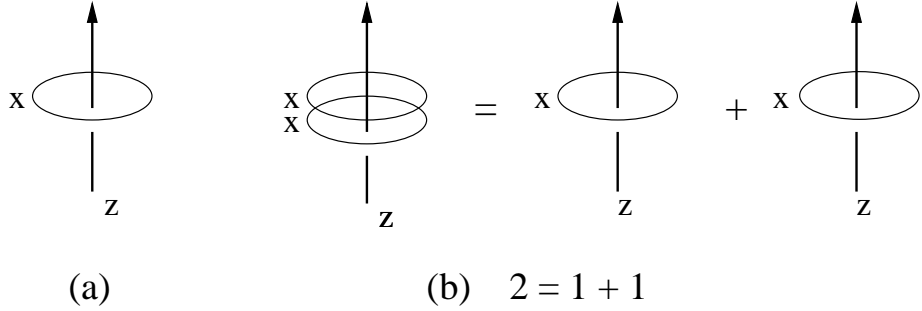}$}
                 }
           }
  \end{equation*}
  \caption{Diagrammatic logic behind Duflo map
  \label{fig_circle_bar}}
 \end{figure}

 \begin{proof}
   Considering the technical complexity, I only want to give the outline of
   the proof in this thesis. For the complete proof, please see \cite{BTT03}.

   First, suppose $H(z;x) \in \mathcal{A}(\ast_z,\ast_x)$ is the disjoint union
   of the diagrams $\Omega_x$ and $\Gamma$. We can show that $H$ is in fact
   the Konsevich integral of the tangle in figure~\eqref{fig_circle_bar}(a).
   Next, we compute the Kontsevich integral of the tangles on both side of
   the diagrammatic equation in figure~\eqref{fig_circle_bar}(b) (so called $2=1+1$),
   and we get
    \begin{equation} \label{oneplusone}
      \Delta_{x_1x_2}H(z;x)=H(z;x_1)\text{\#}_{z} H(z;x_2)
    \end{equation}

   Then for any Jacobi diagrams $D_1, D_2 \in \mathcal{B}$,
   \begin{align} \label{algebra_homo}
     \Phi(D_1 \sqcup D_2) &= \langle H(z;x), (D_1 \sqcup D_2)_{x} \rangle
           \stackrel{~\eqref{duality}}{=}
            \langle \Delta_{x_1x_2}H(z;x), (D_1)_{x_1} \sqcup (D_2)_{x_2} \rangle \notag \\
           &= \langle H(z;x_1)\text{\#}_{z} H(z;x_2),(D_1)_{x_1} \sqcup (D_2)_{x_2} \rangle  \notag \\
           &= \Phi(D_1) \text{\#}_{z} \Phi(D_2)
   \end{align}

  This shows that $\Phi$ is indeed an algebra homomorphism.
\end{proof}

%%%%%%%%%%%%%%%%%%%%%%%%%%%%%%%%%%%%%%%%
%%%%%%% Super Wheeling Theorem %%%%%%%%%
%%%%%%%%%%%%%%%%%%%%%%%%%%%%%%%%%%%%%%%%

 \subsection{Super Version of Wheeling Theorem}
   For a super Lie algebra, the super Jacobi
   identity can be written as:
  \begin{equation}
   (-1)^{degxdegz} [[x,y],z] + (-1)^{degzdegy}[[z,x],y]
   + (-1)^{degydegx} [[y,z],x] =0
  \end{equation}
   or
   \begin{equation}
    [[x,y],z] + (-1)^{degz(degx+degy)}[[z,x],y]
    + (-1)^{degx(degy+degz)} [[y,z],x] =0
   \end{equation}

   Pictorially we can still use figure~\eqref{fig_Jacobi} to
   represent it, plus letting the quadrivalent crossings
   in the diagram represent the sign changes caused by the
   super degree.

\begin{figure}
    \begin{equation*}
      \vcenter{
            \hbox{
                  \mbox{$\includegraphics{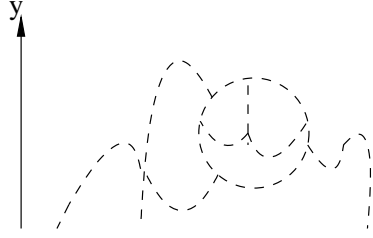}$}
                 }
           }
   \end{equation*}
     \caption{An example of a super diagram}
     \label{fig_super_example}
  \end{figure}

   For any Jacobi diagram $D$ on the x-y
   plane, we can always deforme it to a general position in the upper
   half plane such that all the vertices and crossings have different y-level
   and the tips of legs are placed on the $x$-axis
   (see figure~\eqref{fig_super_example}). Then from
   the bottom to the top, we can decompose the diagram into some forks,
   crossings, cups and caps. If we have a super Lie algebra with an
   invariant metric, we can associate a canonical invariant tensor to
   each of the basic components and then contract their legs to get
   an invariant tensor. It is not hard to show that the result is
   independent of the general position we use for the diagram.

   It is easy to see that the wheeling theorem can be extended to the super
   case without any change.

%%%%%%%%%%%%%%%%%%%%%%%%%%%%%%%%%%%%%%%%%%%%%%
%%%%% Wheeling implies Duflo isomorphism %%%%%
%%%%%%%%%%%%%%%%%%%%%%%%%%%%%%%%%%%%%%%%%%%%%%

\subsection{Wheeling theorem implies the Duflo isomorphism}

\begin{figure}
  \begin{equation*}
   \vcenter{
            \hbox{
                  \mbox{$\includegraphics{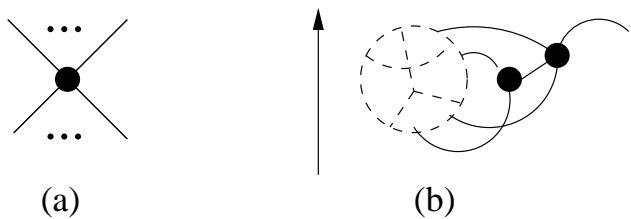}$}
                 }
           }
  \end{equation*}
  \caption{Generalized Graphs that represent arbitrary invariant symmetric tensors
  \label{fig_General_Graph}}
 \end{figure}

\begin{figure}
  \begin{equation*}
   \vcenter{
            \hbox{
                  \mbox{$\includegraphics{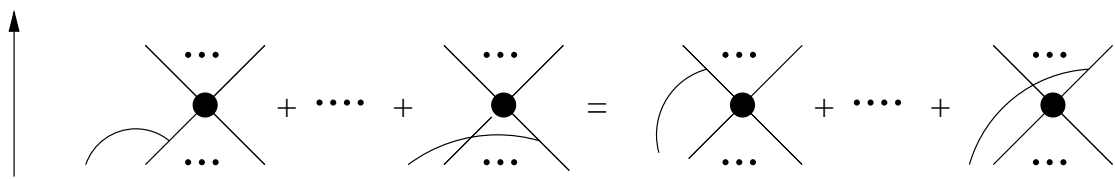}$}
                 }
           }
  \end{equation*}
  \caption{$\g$-invariance of the blobs
  \label{fig_Invariance}}
 \end{figure}

  Although not every invariant symmetric tensor can be represented
  by elements in $\mathcal{B}$, we can introduce some labeled blobs with
  legs to represent arbitrary elements in $S(\g)^{\g}$ as
  in figure~\eqref{fig_General_Graph}(a). So a generalized diagram could
  look like figure~\eqref{fig_General_Graph}(b). The $\g$-invariance of
  the elements represented by those blobs can be drawn as
  diagrammatic relation in figure~\eqref{fig_Invariance}.
  If we put the legs of a generalized diagram on a solid line, it
  will represent an element in $U(\g)^{\g}$.

   The wheeling theorem can be easily extended to this larger
   set of diagrams. So the fact that Duflo isomorphism is an algebra
   isomorphism follows from the diagram below:

      \begin{equation} \label{Super_Duflo}
  \begin{diagram}
   \node{\text{$\mathcal{B}_{blob}$}}
          \arrow[3]{e,tb}{\text{$\Phi$}}{\text{algebra isom.}}
          \arrow{s,l}{\text{surjective}}
   \node[3]{\text{$\mathcal{A}_{blob}$}}
          \arrow{s,r}{\text{surjective}}     \\
   \node{\text{$S(\g)^{\g}$}}
          \arrow[3]{e,tb}{\text{Duflo}}{\text{algebra isom.}}
   \node[3]{\text{$U(\g)^{\g}$}}
  \end{diagram}
 \end{equation} \newline

 %%%%%%%%%%% End of Chapter5
 %Jacobi Diagrams

%%%%%%%%%%%%%%%%%%%%%%%%%%%%%%%%%%%%%%%%%%%%%%%%%%%%%%%%
%%%%% Diagrammatic proof of the quantization map %%%%%%%
%%%%%%%%%%%%%%%%%%%%%%%%%%%%%%%%%%%%%%%%%%%%%%%%%%%%%%%%

\large
\section{\textbf{Diagrammatic analogue of Alekseev-Meinrenken quantization map}}
\normalsize
  In the last chapter, a diagrammatic analogue of Duflo map in
  the world of Jacobi diagrams for any quadratic Lie algebra is shown.
  Although the diagrammatic proof may not
  be easier than the algebraic proof, it provide us a new way to
  handle the invariant tensors which could be useful in some other
  situations.

  The most important evidence I found to believe that
  the quantization map $\mathcal{Q}$
  is equivalent to a super Duflo map is that
  the natural diagrammatic representation of the quantization map
  $\mathcal{Q}$ has the same property as the diagram $H(z,x)$
  for the Duflo map in the preceding chapter.

  \begin{remark}
   The super Lie algebra $\widetilde{T\g[1]}$ for a
   quadratic Lie algebra $\g$ has a natural invariant non-degenerate
    (super)symmetric bilinear form $\widetilde{B}$ defined by:
    $\widetilde{B}(e_a,e_b)= B(e_a,e_b)$,
    $\widetilde{B}(\overline{e_{a}},\overline{e_{b}})=0$,
    $\widetilde{B}(e_a,\overline{e_{b}})=B(e_{a},e_{b})$,
    $\widetilde{B}(e_a,\mathfrak{c})=0$,
    $\widetilde{B}(\overline{e_a},\mathfrak{c})=0$,
    $\widetilde{B}(\mathfrak{c},\mathfrak{c})=1$
    where $\{e_a\}$ and $\{\overline{e_a}\}$ are basis of
    $\widetilde{T\g[1]}^{even}$ and $\widetilde{T\g[1]}^{odd}$
    respectively.
  \end{remark}

   Let's see what diagrammatic analogue of the quantization map $\mathcal{Q}$ should
   be. It is natural to label legs of a Jacobi diagram by
   even (e) or odd (o) with respect to the type of variables they
   represent. We allow partial labelings also. The legs which are
   not labeled by even or odd can be thought of as super legs.

   By the preceding remark, diagrams in $\mathcal{B}$ with legs thus labeled will
   represent a symmetric tensor on $\widetilde{T\g[1]}$ and
   diagrams in $\mathcal{A}$ will represent elements in the
   universal enveloping algebra of $\widetilde{T\g[1]}$. Let
   $\widetilde{\mathcal{B}}$ and $\widetilde{\mathcal{A}}$
   denote the spaces of corresponding Jacobi diagrams
   with legs labeled (or partially labeled) by even and odd. The algebraic
   structures on $\widetilde{\mathcal{B}}$ and
   $\widetilde{\mathcal{A}}$ are the same as $\mathcal{B}$ and
   $\mathcal{A}$ respectively. We can readily call
   them \textit{super labeled Jacobi diagrams}.

   The definition of the diagrammatic differential operator of a
   diagram in this situation is the same as~\eqref{DigramOperator} in
   the previous chapter plus that we only glue even legs to
   even (or super) legs and glue odd legs to odd (or super) legs.

\begin{figure}
  \begin{equation*}
   \vcenter{
            \hbox{
                  \mbox{$\includegraphics{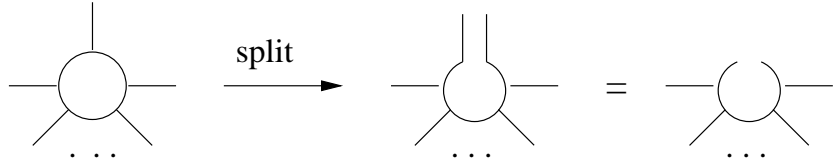}$}
                 }
           }
  \end{equation*}
  \caption{Splitting of a wheel at a leg
  \label{figure_splitting}}
 \end{figure}

   Next, we want to use the Jacobi diagram of $j^{\half}$ to
   define the diagram for the tensor $\exp(\half T_{ab}(x)\iota_{a}\iota_{b})$.
   Since $T_{ab}(x)= (\ln(j)'(ad_{x}))_{ab}$, we need to define
   the \textit{derivative of a Jacobi diagram}.
   First, we call the way of changing a diagram in figure~\eqref{figure_splitting}
   \textit{splitting a wheel at a leg}.
   Since the derivative of $x^n$ is $nx^{n-1}$,
   so it is natural to define the diagrammatic derivative $\mathcal{P}(C)$ of a
   wheel $C$ to be the diagrammatic sum of all possible ways of splitting $C$ at
   its legs. For the disjoint union of two wheels
   $C_1$ and $C_2$, we define $\mathcal{P}(C_1 \sqcup C_2) =
   \mathcal{P}(C_1) \sqcup \mathcal{P}(C_2)$. With this definition, we have
   $\mathcal{P}(\exp^{\sum_{i}C_{i}})=\exp^{\sum_{i}\mathcal{P}(C_{i})}$ for
   a collection of wheels $C_i$. In addition, we define the action of
   $\mathcal{P}$ on a connected diagram other than wheels to be
   trivial (i.e. kill it).

\begin{figure}
  \begin{equation*}
   \vcenter{
            \hbox{
                  \mbox{$\includegraphics{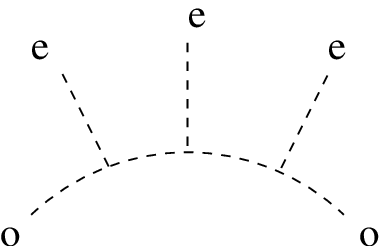}$}
                 }
           }
  \end{equation*}
  \caption{$ad^{n}(x)\iota_{a}\iota_{b}$, $n=3$
  \label{figure_T_ab}}
 \end{figure}

\begin{figure}
  \begin{equation*}
   \vcenter{
            \hbox{
                  \mbox{$\includegraphics{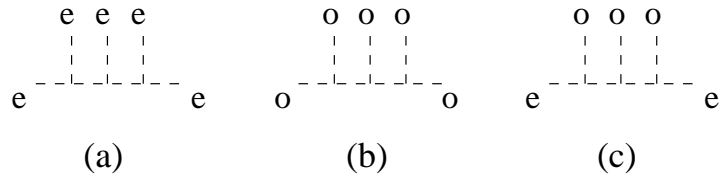}$}
                 }
           }
  \end{equation*}
  \caption{Three other possible cases of splitting
  \label{figure_3cases}}
 \end{figure}

\begin{figure}
  \begin{equation*}
   \vcenter{
            \hbox{
                  \mbox{$\includegraphics{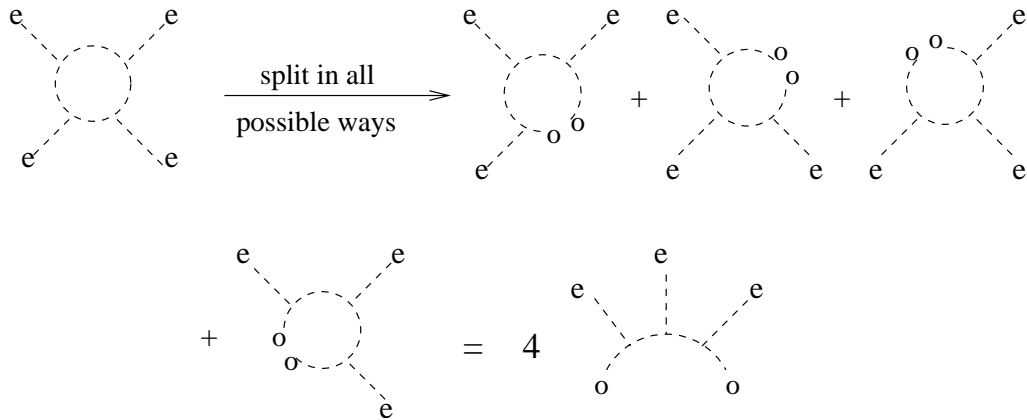}$}
                 }
           }
  \end{equation*}
  \caption{Diagrammatic derivative of a wheel
  \label{figure_Split}}
 \end{figure}

  The diagram for the Duflo map of $\widetilde{T\g[1]}$ has super
  legs. A super leg consist of an even leg and an odd leg.
  Then there are four possible cases of splitting.
  The first case is shown in figure~\eqref{figure_T_ab} which
  corresponds to $ad^{n}(x)\iota_{a}\iota_{b}$. The other three
  cases (see figure~\eqref{figure_3cases}) are in fact all zero algebraically.
  The case(a) vanishes because $T_{ab}$ is skew-symmetric,
  $T_{ab}\pardev{}{\mu^{a}}\pardev{}{\mu^{b}}= -T_{ba}\pardev{}{\mu^{a}}\pardev{}{\mu^{b}}$,
  and the cases (b) and (c) vanish obviously. The diagrammatic
  derivative of a $4$-leg wheel is shown in figure~\eqref{figure_Split}.

   \begin{remark}
   In figure~\eqref{figure_T_ab}, the even legs in fact correspond to
   $v^a(\half f_{abc}\bar{e_{b}}\bar{e_{c}})$ after the super variable change.
   \end{remark}

   By the Taylor expansion of $\ln(j^{\half})$ in~\eqref{eq:Bernoulli}, we get the
   diagrammatic representation $\Psi_{x}$ of $\exp(\half T_{ab}(x)\iota_{a}\iota_{b})$.
  \begin{align}
    \Psi_{x} &= \exp(\half T_{ab}(x)\iota_{a}\iota_{b}) \notag \\
             &= \exp \left(
                  \frac{1}{24}\: \includegraphics[scale=0.8]{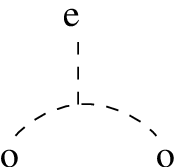}
              - \frac{1}{1440}\: \includegraphics[scale=0.8]{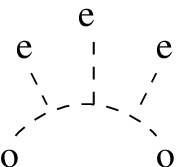}
                       + \ldots
                     \right)
  \end{align}

   In addition, the map $\chi \otimes q$ in quantization map
   $\mathcal{Q}$ can still be represented by $\Gamma$ (see
   figure~\eqref{fig_Chi}). The legs attached to the z-line are
   not labeled by even or odd because we allow them to be connected to both even and odd
   type of legs. Notice we won't have any legs directly attached to the z-line after the
   contraction $\partial_{\Gamma}$ with a Jacobi diagram
   (see the definition of $\partial_{\Gamma}$ in~\eqref{Diagram_Drivative}).

   Let $\widetilde{H}(z;x)$ to be the disjoint union of $\Gamma ,\Omega_{x}$
   and $\Psi_{x}$(labeling ignored). by the definition of
   $\mathcal{P}$, we have the following lemma.
   \begin{lemma}
    $\widetilde{H}(z;x)=(I+\mathcal{P})H(z;x)$, where $I$ is the
    identity map.
   \end{lemma}

   Notice
   \textit{labeling} and \textit{coloring} are different actions on diagrams.
   More precisely, labeling a leg by even or odd and coloring a leg by x or y are
   independent to each other.
   Observe that: the splitting operation $\mathcal{P}$ on wheels commutes
   with the operation $\Delta_{xy}$ which is just ignore the original coloring of a
   diagram and coloring all legs of a diagram by two colors $x,y$ in all
   possible ways (see~\eqref{coloring} for definition).
   Then we have the following:
   \begin{align}
     \Delta_{x_1x_2}\widetilde{H}(z;x)
     &= \Delta_{x_1x_2} (I+\mathcal{P})H(z;x)
      = \Delta_{x_1x_2} (I+\mathcal{P}) \Delta_{x_1x_2} H(z;x) \notag \\
     &= \Delta_{x_1x_2}(I+\mathcal{P}) \{ H(z;x_1) \text{\#}_{z} H(z;x_2)
        \}  \notag \\
     & = \{ \Delta_{x_1x_2}  (I+\mathcal{P})H(z;x_1) \} \text{\#}_{z}
         \{ \Delta_{x_1x_2}  (I+\mathcal{P})H(z;x_2) \}  \notag \\
     &= \widetilde{H}(z;x_1) \text{\#}_{z} \widetilde{H}(z;x_2)
   \end{align}

   The second equality needs a little explanation. Notice the
   action of $\mathcal{P}$ will introduce some new odd legs which
   are not colored by x or y yet, so we need to remove the color we just did and
   do the coloring again to guarantee each leg is colored.

   Next, using the same argument in~\eqref{algebra_homo}, we can
   easily prove the following theorem which can be thought of as the
   diagrammatic analogue of quantization map $\mathcal{Q}$ for Weil algebras.

   \begin{thm}
   If we label the legs of $\Omega_x$ by even,
   the map $\widetilde{\Phi} = \partial_{\Gamma} \circ \partial_{\Omega_x}
    \circ \partial_{\Psi_{x}}  : \widetilde{\mathcal{B}}
   \longrightarrow \widetilde{\mathcal{A}}$
   is an algebra homomorphism with respect to
   the algebraic structures of $\widetilde{\mathcal{B}}$ and
   $\widetilde{\mathcal{A}}$. \newline
   \end{thm}

%%%%% End of chapter 6 %%%%%%%%%
 %Diagrammtic analogue of Quantization map

\bibliographystyle{abbrv}
\bibliography{thesis}

\end{document}